\documentclass[a4paper]{amsart}
\usepackage{amsmath,latexsym}
\usepackage[mathscr]{eucal}
\usepackage{amssymb}

\newcommand{\gs}{\mathfrak{S}}
\newtheorem*{HypoSY}{Hypothesis $S(Y)$}
\newtheorem{theorem}{Theorem}
\newtheorem{lemma}{Lemma}

\allowdisplaybreaks
\theoremstyle{remark}

\newtheorem*{rema}{Remark}

\begin{document}

\normalsize

\renewcommand{\thefootnote}{}

\title{Primes in Tuples II}

\author{D. A. Goldston, J. Pintz and
C. Y. Y{\i}ld{\i}r{\i}m}

\date{\today}

\numberwithin{equation}{section}

\begin{abstract}
We prove that 
$$
\liminf\limits_{n \to \infty}
\frac{p_{n+1} - p_n}{\sqrt{\log p_n}(\log\log p_n)^2}
\ < \infty,
$$
where $p_n$ denotes the $n$\textsuperscript{th} prime. Since on average $p_{n+1} - p_n$ is asymptotically $\log p_n$, this shows that we can always find pairs of primes much closer together than the average. We actually prove a more general result concerning the set of values taken on by the differences $p-p'$ between primes which includes the small gap result above.
\end{abstract}

\maketitle
\footnotetext{The first author was supported in part by an NSF Grant,
the second author by OTKA grants No. K~67676, T~43623, T~49693
and the Balaton program,
the third author by T\"UBITAK.}

\section{Introduction}
\label{sec:1}

In the first paper in this series \cite{GPY} we proved that, letting $p_n$ denote the
$n$\textsuperscript{th} prime,
\begin{equation}
\Delta = \liminf_{n\to \infty} \frac{p_{n+1} - p_n}{\log p_n} =0 ,
\label{eq:1.1}
\end{equation}
culminating 80 years of work on this problem. Since the average spacing $p_{n+1} - p_n$ in the sequence of primes is asymptotically $\log p_n$, this result showed for the first time that the prime numbers do not eventually become isolated from each other in the  sense that there will always be pairs of primes closer than any  fraction of the average spacing. 
For the history of this problem, we refer the reader to \cite{GPY} and \cite{Pi4}. 

The information about primes used to obtain \eqref{eq:1.1} is contained in the Bombieri-Vinogradov theorem. Let
\begin{equation}\theta(N;q,a) = \sum_{\substack{n\le N \\ n\equiv a (\text{mod}\,
q)}}\theta(n), \quad \text{where} \ \
\theta(n) = \left\{
\begin{array}{ll}
  \log n, &\text{if } n \text{ is  prime},\\
  0, &\text{otherwise}.
\end{array}
\right. \label{eq:1.2}
\end{equation}
 The Bombieri-Vinogradov theorem states that for any  $A>0$ there is a $B=B(A)$ such that, for $Q=N^{\frac{1}{2}}(\log N)^{-B}$,
\begin{equation}
\sum_{q\le Q}  \max_{\substack{a\\  (a,q)=1}}  \left| \theta(N;q,a)-   \frac{N}{\phi(q)}  \right| \ll \frac{N}{(\log N)^A}. \label{eq:1.3}
\end{equation}
Thus the primes tend to be equally distributed among the arithmetic progressions modulo $q$ that allow primes, and this holds for the primes up to $N$ and at least for almost all the progressions with modulus $q$ up to nearly $N^{\frac12}$. The principle can be quantified by saying that the primes have an \emph{admissible level of distribution} $\vartheta$ (or \emph{satisfy a level of distribution} $\vartheta$) if \eqref{eq:1.3} holds for any $A>0$ and any $\epsilon >0$ with
\begin{equation} Q=N^{\vartheta -\epsilon}.  \label{eq:1.4} \end{equation}  
Elliott and Halberstam \cite{EH} conjectured that the primes have the maximal admissible level of distribution 1, while by the Bombieri-Vinogradov theorem we have immediately that $1/2$ is an admissible level of distribution for the primes. In \cite{GPY} we proved that \emph{if} the primes satisfy a level of distribution $\vartheta > \frac12$ then there is an absolute constant $M(\vartheta)$ for which
\begin{equation}  p_{n+1}-p_{n} \le  M(\vartheta),  \quad \text{ for infinitely many $n$.} \label{eq:1.5}\end{equation}
In particular assuming the Elliott-Halberstam conjecture (or just $\vartheta \ge 0.98$) then
\begin{equation}  p_{n+1}-p_{n} \le  16,  \quad \text{ for infinitely many $n$.} \label{eq:1.6}\end{equation}
These are surprising results because they show that going beyond an admissible level of distribution 1/2 implies there are infinitely often bounded gaps between primes, and therefore questions as hard as the twin prime conjecture can nearly be dealt with using this type of information. 

Since we obtained our results in 2005 there has been no further progress toward \eqref{eq:1.5}, and it appears now that an extension of the Bombieri-Vinogradov theorem of sufficient strength to obtain bounded gaps between primes will require some basic new ideas. One can also pursue improving the approximations we used or improving the method used to detect primes, and again this now appears to require some essentially new idea. Our goal in this paper is to extend the current method as much as possible in order to obtain strong quantitative results. In particular we obtain the following quantitative version of \eqref{eq:1.1}.

\begin{theorem}
\label{th:A}
The differences of consecutive primes satisfy
\begin{equation}
\liminf_{n \to\infty} \frac{p_{n+1} - p_n}{\sqrt{\log
p_n}(\log\log p_n)^2} < \infty .
\label{eq:1.7}
\end{equation}
\end{theorem}
This result is remarkable in that it shows that there exist pairs of primes nearly within the square root of the average spacing. By comparison, the best result for large gaps between primes \cite{Pi3} is that 
\begin{equation}
\limsup_{n \to\infty} \frac{p_{n+1} - p_n}{(\log
p_n)(\log\log p_n) (\log\log\log p_n)^{-2}(\log\log\log\log p_n)}
\ge 2e^\gamma ,
\label{eq:1.8}
\end{equation}
where $\gamma$ is Euler's constant.
Thus the best large gap result produces gaps larger than the average by a factor a bit smaller than $\log\log p_n$, while now the small gaps are smaller than the average by a factor a bit bigger than $(\log p_n)^{-\frac12}$. In this sense the small gap result now greatly surpasses the large gap result. 
(There are conflicting conjectures on how large the gap between consecutive primes can get, but all of these conjectures suggest that there can be gaps at least as large as $c (\log p_n)^2$ for some constant $c$.) 

\bigskip

This paper is organized as follows. In Section 2 we present a generalization of Theorem 1 which applies to many interesting situations and which is the result we will prove in this paper. Unlike in \cite{GPY}, to obtain our results we need to take into account the possibility of exceptional characters associated with  Landau-Siegel zeros. In Theorem 2 we assume that there are no Landau-Siegel zeros in a certain range and are able to obtain our main results including the gaps between primes in Theorem 1 in intervals $[N,2N]$ for all sufficiently large $N$. Next, using the Landau-Page Theorem we can find a sequence of ranges which avoid possible Landau-Siegel zeros. Thus we obtain Theorem 3 which unconditionally gives the same results as Theorem 2 but without being able to localize them to a dyadic interval.   The proof of these theorems requires  substantial refinements of the methods of \cite{GPY}, and in Section 3 we will discuss some of these refinements and how they arise. The main technical tools needed in our proof, Theorems 4 and 5, are stated in Section 4. The proof of Theorems 4 and 5 take up Sections 5 through 13. In proving Theorem 5 in our general setting we need a modified Bombieri-Vinogradov theorem which is the topic of Section 12.  Our method, as in \cite{GPY}, requires a result on the average of the singular series. In \cite{GPY} the well-known result of Gallagher \cite{Ga} was used, but in our current setting this result is not applicable, and therefore in Section 14 we prove a new result well adapted for our needs. With Theorems 4 and 5 in hand together with the new singular series average result, the proof of Theorems 2 and 3 is completed in Section 15.

\subsection*{Notation}
In the following $c$ and $C$ will denote (sufficiently) small
and (sufficiently) large absolute positive constants,
respectively, which have been chosen appropriately.
This is also true for constants formed from $c$ or $C$ with
subscripts or accents.
We will allow these constants to be different at different occurences.
Constants implied by pure $o$, $O$, $\ll$ symbols will be
absolute, unless otherwise stated.
The $\nu$ times iterated logarithm will be denoted by
$\log_\nu N$. $\mathcal{P}$ denotes the set of primes.

\section{A generalization of Theorem~1}
\label{sec:uj2}

Our method will allow us to prove a generalization of
Theorem~1 where instead of seeking two neighboring primes of the
form $n + i$, $n + j$ with
\begin{equation}
1 \leq i < j \leq h, \qquad
h = h(n) = C\sqrt{\log n} (\log\log n)^2,
\label{eq:uj2.1}
\end{equation}
we look for two primes of the form $n + a_i$, $n + a_j$ where
\eqref{eq:uj2.1} is satisfied and
\begin{equation}
\mathcal A = \{ a_i\}^h_{i = 1} \subset [1, n]
\label{eq:uj2.2}
\end{equation}
is an arbitrarily given set of integers.

We remark that an extension of this
type is a trivial consequence of the prime number theorem if
\begin{equation}
h' = h'(n) > (1 + c) \log n, \quad c > 0 \text{ fixed},
\label{eq:uj2.3}
\end{equation}
but that none of the earlier methods of Erd\H os \cite{Erdos}, Bombieri--Davenport
\cite{BD} and Maier \cite{Ma} which produce small gaps between primes
seem capable of proving a result of this type for any
function satisfying
\begin{equation}
h '' = h''(n) < (1 - c) \log n, \quad c > 0 \text{ fixed.}
\label{eq:uj2.4}
\end{equation}

According to a conjecture of de Polignac \cite{Poli1849b}
from 1849,  every even number may be
written as the difference of two primes.
Although we know that this is true for almost all even numbers,
there is no known way to specify these values.
Our generalization makes a first step in this
direction by proving that we can explicitly find sparse
sequences $\mathcal{A}$ such that infinitely many of the elements
\begin{equation}
\mathcal A - \mathcal A
\label{eq:uj2.2.4}
\end{equation}
are 
differences of two primes, i.e.
\begin{equation}
|(\mathcal P - \mathcal P) \cap (\mathcal A - \mathcal A)| = \infty.
\label{eq:uj2.5}
\end{equation}
(Here we make use of the usual notation that for sets $\mathcal A$ and $\mathcal B$, 
$ \mathcal A -\mathcal B = \{ a -b : a \in \mathcal A, b \in \mathcal B\}.$)
Some sequences for which our method applies are:
\begin{alignat}{2}
\mathcal A &= \{k^m\}^\infty_{m = 1}, \ &  k &\geq 2 \text{
fixed } (k \in \mathbb N), \label{eq:uj2.6}\\
\mathcal A &= \{k^{x^2 + y^2}\}^\infty_{x,y = 1}, \ &  k &\geq 2 \text{
fixed } (k \in \mathbb N), \label{eq:uj2.7}\\
\mathcal A &= \{k^{f(x,y)}\}^\infty_{r = 1}, \ &  k &\geq 2 \text{
fixed } (k \in \mathbb N), \label{eq:uj2.8}
\end{alignat}
where the value set $\mathcal R = \{r \in \mathbb N;\ \exists x,y : f(x,y) =
r\}$ satisfies
\begin{equation}
\mathcal R(X) = |\{ m;\ m \leq X,\ m \in \mathcal R\}| > C'
\sqrt{X} \log^2 X
\label{eq:uj2.9}
\end{equation}
(this would happen e.g.\ for $f(x,y) = x^2 + y^m$ with arbitrary $m \geq 2$, 
and for $f(x,y)=x^3 + y^3$), or in general any set of type
\begin{equation}
\mathcal A = \{k^{r_j}\}^\infty_{j = 1},\quad k \geq 2
\text{ fixed } (k \in \mathbb N),
\label{eq:uj2.10}
\end{equation}
if $\mathcal R = \{r_j\}^\infty_{j = 1} \subseteq \mathbb N$
satisfies the density condition \eqref{eq:uj2.9}.
Among these sets the only trivial one is $2^m$, that is
\eqref{eq:uj2.6} for $k = 2$ which has $[\log X / \log 2]$
elements below $X$, thereby corresponding to the
case~\eqref{eq:uj2.3}.

Unfortunately, the possible existence of Landau-Siegel zeros (cf.\
Section~12) makes it impossible to formulate a localized version of our 
result for all $n \in \mathbb N$ satisfying the conditions
\eqref{eq:uj2.1} and \eqref{eq:uj2.2}.
Even in the special case when $\mathcal A$ is an interval
we cannot guarantee the existence of
gaps of size $O\bigl(\sqrt{\log N} \log^2_2 N\bigr)$
between primes in any interval of type $[N, 2N]$ for any large~$N$.

The first formulation of our main result assumes there are no exceptional characters in a certain range. (This is Hypothesis
$S(Y)$ from Section~12, with $Y = Y(N) = \exp\bigl(3\sqrt{\log
N}\bigr)$.)

\begin{theorem}
\label{th:uj2}
Let us suppose that an $N > N_0$ is given such that for any real
primitive character $\chi\mod q$, $q \leq \exp\bigl(3 \sqrt{\log
N}\bigr)$ we have
\begin{equation}
L(s, \chi ) \neq 0 \ \text{ for } \ s \in \left(1 -
\frac1{9 \sqrt{\log N}}, 1 \right] .
\label{eq:uj2.11}
\end{equation}

Let $\mathcal A = \mathcal A_N = \{a_i\}^h_{i = 1} \subseteq [1,
N] \cap \mathbb N$ be arbitrary $(a_i \neq a_j)$ with
\begin{equation}
h \geq C \sqrt{\log N} \log^2_2 N,
\label{eq:uj2.12}
\end{equation}
where $C$ is an appropriate absolute constant.
Then there exists $n \in [N, 2N]$ such that at least two numbers
of the form
\begin{equation}
n + a_i, \ n + a_j, \qquad (1 \leq i < j \leq h),
\label{eq:uj2.13}
\end{equation}
are primes.
\end{theorem}

In Section~12 we show  that 
\eqref{eq:uj2.11}, i.e. Conjecture $S(Y(N))$, is true for an infinite sequence $N = N_\nu
\to \infty$, and thus  Theorem~\ref{th:uj2} 
implies Theorem~\ref{th:A} by choosing $N = N_\nu$ and
\begin{equation}
\mathcal A = \mathcal A_N = \{1,2,\dots, h\}
\label{eq:uj2.13}
\end{equation}
with $h = \left\lceil C \sqrt{\log N} \log^2_2 N\right\rceil$.
A more general formulation of this result which covers the
special cases mentioned in \eqref{eq:uj2.6}--\eqref{eq:uj2.8} and \eqref{eq:uj2.10} as well as
Theorem~\ref{th:A} can be stated as follows.

\begin{theorem}
\label{th:uj3}
Let $\mathcal A \subseteq \mathbb N$ be an arbitrary sequence satisfying
\begin{equation}
\mathcal A(N) = \bigl|\{n; n \leq N, n \in \mathcal A\}\bigr| > C \sqrt{
\log N} \log^2_2 N \text{ for } N > N_0.
\label{eq:uj2.14}
\end{equation}
Then infinitely many elements of $\mathcal A - \mathcal
A$ can be written as the difference of two primes, that is,
\begin{equation}
|(\mathcal P - \mathcal P) \cap (\mathcal A - \mathcal A)| = \infty.
\label{eq:uj2.15}
\end{equation}
\end{theorem}

\section{Some Initial Considerations}
\label{sec:2}

The main tool of our method is an approximation for prime tuples and almost prime tuples. Consider
the tuple $(n+h_1,n+h_2, \ldots, n+h_K)$ as $n$ runs over the integers. If these
$K$ values are all primes for some $n$ then we call this a prime tuple, and we wish to examine the existence of prime tuples.  A first consideration is that the set of shifts
\begin{equation}
\mathcal{H}=\{h_1,h_2, \ldots, h_K\}, \quad
\textrm{with }  \  h_i \neq h_j \
(\text{if } i \neq j),
\label{eq:2.1}
\end{equation}
imposes divisibility conditions on the components of the tuple which can effect the likelihood of obtaining prime tuples or even preclude the possibility of more than a single prime tuple. Specifically, let
 $\nu_p(\mathcal{H})$ denote the number of distinct
residue classes modulo $p$ occupied by the elements of
$\mathcal{H}$, and for squarefree integers $d$  extend this
definition to $\nu_d(\mathcal{H})$ multiplicatively. 
The singular series for the set $\mathcal{H}$ is defined to be
\begin{equation}
\mathfrak{S}(\mathcal{H}) = \prod_p\left(
1-\frac{1}{p}\right)^{-K}\left(1 -
\frac{\nu_p(\mathcal{H})}{p}\right) .
\label{eq:2.2}
\end{equation}
If $\mathfrak{S}(\mathcal{H})\neq 0$ then $\mathcal{H}$ is called
\emph{admissible}. Thus $\mathcal{H}$ is admissible if and only
if $\nu_p(\mathcal{H})<p$ for all $p$, while if $\nu_p(\mathcal{H})=p$ then one component of the tuple is always divisible by $p$ and there can be at most one prime tuple of this form. 
Hardy and Littlewood \cite{HL} conjectured an asymptotic formula for the number of prime tuples $(n+h_1,n+h_2, \ldots , n+h_K)$, with
$1\le n\le N$, as $N\to \infty$. Letting
\begin{equation}
\theta(n) = \begin{cases}
\log n, &\text{if $n$ is prime,}\\
0, &\text{otherwise;}
\end{cases}
\label{eq:2.23}
\end{equation} we define
\begin{equation}
\Lambda(n;\mathcal{H}) := \theta(n+h_1)\theta(n+h_2)\cdots
\theta(n+h_K) \label{3.3}
\end{equation}
and use this function to detect prime tuples.
The Hardy--Littlewood prime-tuple conjecture is the asymptotic formula
\begin{equation}
\sum_{n\le N}\Lambda(n;\mathcal{H}) = N (\gs(\mathcal{H})+o(1)),
\quad  \mbox{as \ $N\to \infty$,}\label{3.5}
\end{equation}
which is trivial if $\mathcal{H}$ is not admissible, but is otherwise only known to be true in the case $K=1$ which is the prime number theorem.

The starting point for our method in \cite{GPY} is to find approximations of $\Lambda(n;\mathcal{H})$ for which we can obtain asymptotic formulas similar to \eqref{3.5}. A further essential idea is that rather than approximating just prime tuples, we should approximate almost-prime $K$-tuples with a total of $\le K+\ell$ prime factors in all the components, which if $0\le \ell \le K-2$ guarantees at least two of the components are prime.
The almost prime tuple approximation used in
\cite{GPY} and which we also use here is  
\begin{equation}
\Lambda_R(n;\mathcal H, \ell) : = \frac{1}{(K + \ell)!} \sum_{\substack{d
\mid P_{\mathcal H}(n)\\ d \leq R}} \mu(d) \left(\log
\frac{R}{d}\right)^{K + \ell},
\label{eq:2.3}
\end{equation}
where  $|\mathcal H|=K$, and
\begin{equation}
\mathcal P_{\mathcal H}(n) := (n + h_1)(n + h_2) \dots (n + h_K).
\label{eq:2.4}
\end{equation}
Our method for proving \eqref{eq:1.1} in \cite{GPY} is based on a comparison of the two sums
\begin{equation} \sum_{n\leq N} \Lambda_{R} (n; \mathcal{H} , \ell)^2 \quad \text{and} \quad \sum_{n\leq N} \theta(n+h_0)\Lambda_{R} (n; \mathcal{H} , \ell)^2.\label{3.8new}\end{equation}

An asymptotic formula for the first sum can be obtained if $R\le N^{1/2-\epsilon}$, while for the second sum  we can use an admissible level of distribution of primes $\vartheta$ to obtain an asymptotic formula when $R\le N^{\vartheta/2-\epsilon}$. In \cite{GPY} it was assumed that $K$ and $\ell$ are fixed, i.e. independent of $N$. Using these asymptotic formulas we can now evaluate
\begin{equation}
\mathcal{S}_R
:= \sum_{n=N+1}^{2N}
\left( \sum_{1\le h_0\le h} \theta (n+h_0) -  \log 3N \right)
\sum_{\substack{1\le h_1,h_2, \ldots ,h_K\le h \\
\text{distinct}}}
\Lambda_{R}(n;\mathcal{H}, \ell)^2 ,
\label{3.8}
\end{equation}
If $\mathcal{S}_R>0$ then the sum over $h_0$ must have at least two non-zero terms and thus there must be some $n$ and $h_i\neq h_j$ such that $n+h_i$ and $n+h_j$ are both prime. We find with $\vartheta = 1/2$ and $h = \lambda \log N$ with any fixed $\lambda >0$ that we can choose $K$ and $\ell$ for which $\mathcal{S}_R>0$, which proves \eqref{eq:1.1}. In order to obtain this for any arbitrarily small $\lambda >0$, the fixed $K$ and $\ell$ are chosen  sufficiently large in an appropriate way. 

To obtain quantitative bounds to replace \eqref{eq:1.1}, the first step is to obtain asymptotic formulas which are uniform in $K$ and $\ell$ so that these can be chosen as  functions of $N$ that go to infinity with $N$. One also needs explicit error terms, and these error terms arise not only from lower order terms and prime number theorem type error terms, but also in \eqref{3.8new} from the Bombieri-Vinogradov theorem error terms. 

We will now establish the relations between our parameters that will be used throughout the paper. Recalling the set $\mathcal{A}$ from \eqref{eq:uj2.2}, we will always take  
$\mathcal H\subset \mathcal A$. Next,  $R$ and $\ell$ will be chosen as
\begin{equation}  K \leq h, \quad
\ell \asymp \sqrt K, \quad
R := (3N)^\Theta = (3N)^{1/4 - \xi}, \qquad \xi = o(1).
\label{eq:2.5}
\end{equation}
We will make use of two important parameters $U$ and $V$ defined by
\begin{equation}
V := \sqrt{\log N}, \ U = e^V,
\label{eq:2.13}
\end{equation}
and will choose $K $ later to be slightly smaller than $V$. 
We next denote the product of primes not exceeding $V$ by 
\begin{equation}
P := \prod_{p \leq V} p,
\label{eq:2.14}
\end{equation}
where $p$ will always denote primes.

As just mentioned above, our present treatment requires a much more delicate analysis of the error terms than in
\cite{GPY}, and therefore we make an initial simplification to facilitate this analysis. In \cite{GPY} the irregular
behavior of $\nu_p(\mathcal
H)$ for small  primes greatly complicated the estimate of the function
$G(s_1, s_2)$ and its partial derivatives. We can avoid these difficulties,
 at least for primes
dividing $P$, by proceeding somewhat similarly to Heath-Brown in
 \cite{HB2}. 
We  call a residue class $a(\text{\rm mod}\,
P)$
\emph{regular with respect to $\mathcal H$  and $P$} if
\begin{equation}
(P, P_{\mathcal H}(a)) = 1
\label{eq:2.15}
\end{equation}
and denote by $A(\mathcal H) = A_P(\mathcal H)$ the set of all
regular residue classes $\text{\rm mod}\, P$. Thus
\begin{equation}
A(\mathcal H) := \big\{ a; \ 1 \leq a \leq P;\ (P, P_{\mathcal
H}(a)) = 1\big\}.
\label{eq:2.16}
\end{equation}
The number of regular residue classes $\text{\rm mod}\, P$ is clearly
\begin{equation}
|A(\mathcal H)| = \prod_{p \mid P} (p - \nu_p(\mathcal H))
\label{eq:2.17}
\end{equation}and their proportion of all the residue classes $\text{\rm mod}\, P$ is
\begin{equation}
\frac{|A(\mathcal H)|}{P} = \prod_{p\mid P} \left(1 -
\frac{\nu_p(\mathcal H)}{p}\right),
\label{eq:2.18}
\end{equation}
which is positive if $\mathcal H$ is admissible.
Thus in particular for a given $\mathcal H$ and all $P$ there exists at least
one regular residue class $\text{\rm mod}\,P$, if and only if $\mathcal H$ is admissible.

With this notation, we now consider the sums
\begin{equation}
\sum^{2N}_{\substack{n = N + 1\\ n \in A(\mathcal H_1) \cap A(\mathcal
H_2)}} \Lambda_R(n; \mathcal H_1, \ell) \Lambda_R(n; \mathcal
H_2, \ell)
\label{eq:2.21}
\end{equation}
and
\begin{equation}
\sum^{2N}_{\substack{n = N + 1\\ n \in A(\mathcal H_1) \cap A(\mathcal
H_2)}} \Lambda_R(n; \mathcal H_1, \ell) \Lambda_R(n; \mathcal
H_2, \ell) \theta(n + h_0)
\label{eq:2.22}
\end{equation}
with $h_0 \in [1, h]$, which are asymptotically evaluated in Theorems~\ref{ujth:1} and \ref{ujth:2},
respectively. A new feature in the proof of these theorems which does not occur in \cite{GPY} is that
 \eqref{eq:2.21} and \eqref{eq:2.22} are first evaluated for
each residue class
$a(\text{\rm mod}\, P)$ with
\begin{equation}
a \in A(\mathcal H_1) \cap A(\mathcal H_2) = A(\mathcal H_1 \cup
\mathcal H_2)
\label{eq:2.24}
\end{equation}
separately, and then the results are added over all regular residue
classes modulo~$P$. It turns out that the asymptotic main term (and even
secondary terms) are independent of the particular choice of
the regular residue class $a$, so this summing presents no difficulty.
However, the restriction of the values of $n$ to a
single residue class $a(\text{\rm mod}\, P)$ in \eqref{eq:2.22} requires a 
stronger form of the Bombieri--Vinogradov theorem (cf.\ Section~12).

To detect primes, in place of  \eqref{3.8} we consider
\begin{equation}
S'_R(N, K, \ell, P) := \frac{1}{Nh^{2K + 1}} \sum^{2N}_{n = N +
1} \bigg(\sum_{\substack{p\\
p - n \in {\mathcal A}}} \log p - \log 3N \bigg)
\big(\Psi'_R (K, \ell, n, h) \big)^2,
\label{eq:2.19}
\end{equation}
where
\begin{equation}
\Psi'_R(K, \ell, n, h) := \sum_{\substack{\mathcal H, |\mathcal H| =
K\\ n \in A(\mathcal H)}} \Lambda_R(n; \mathcal H, \ell).
\label{eq:2.20}
\end{equation}

On applying Theorems 4 and 5 we can asymptotically evaluate $S'_R$, which we carry out in Section 15. One condition that arises from the main terms is that in order to prove the existence of prime pairs in
intervals of length $h$ we need 
\begin{equation}
h > \frac{C \log N}{K}.
\label{eq:2.6}
\end{equation}
Since $K \leq h$ this immediately implies
that
\begin{equation}
h > C \sqrt{\log N}.
\label{eq:2.10}
\end{equation}
Our goal is to take $h$ as small as possible, and therefore we can not obtain anything better than \eqref{eq:2.10} when using the approximation in \eqref{eq:2.3} together with \eqref{eq:2.19}.  Apart from powers of  $\log_2 N$, we are able to prove our results for $h$ of this size.

Our actual choices for $K$ and $h$ are
\begin{equation}
K \leq c_1 \frac{\sqrt{\log N}}{\log^2_2 N}, \quad
h = \frac{25 \log N}{K} \geq \frac{25}{c_1} \sqrt{\log N}
\log^2_2 N,
\label{eq:2.12}
\end{equation}
with a sufficiently small explicitly calculable absolute
constant $c_1$ (to be chosen later).
We will need the error terms in Theorems 4 and 5 to be
uniform in $K$ with a relative error of size $\eta_1$ satisfying
\begin{equation}
\eta_1 < \frac{c}{\sqrt K}.
\label{eq:2.8}
\end{equation}
However, we do not achieve this for all admissible pairs $\mathcal H_1$ and $\mathcal
H_2$ of size $K$. Instead, for 
 all admissible pairs $\mathcal H_1$,
$\mathcal H_2$ we obtain a weaker error term, but
if
\begin{equation}
K - |\mathcal H_1 \cap \mathcal H_2| \ll \sqrt{K}.
\label{eq:2.25}
\end{equation}
then we do obtain the error estimate in \eqref{eq:2.8}.
This turns out to be sufficient for our proof, since such pairs $\mathcal H_1,
\mathcal H_2$ will be dominant in~\eqref{eq:2.19}.


\section{Two basic theorems}
\label{sec:3}

In the following let $N$ be a sufficiently large integer, $c_1$
a sufficiently small positive constant,
\begin{equation}
K \leq c_1 \frac{\sqrt{\log N}}{(\log_2 N)^2},
\label{eq:3.1}
\end{equation}
\begin{equation}
K \ll k_1, k_2 \leq K, \quad \sqrt K \ll \ell_1, \ell_2 \ll
\sqrt K.
\label{eq:3.2}
\end{equation}
We will consider sets
$\mathcal H := \mathcal H_1 \cup \mathcal H_2$,
$\mathcal
H_1, \mathcal H_2 \subseteq [1, N]$ of sizes
\begin{equation}
|\mathcal H_i| = k_i, \quad |\mathcal H_1 \cap \mathcal H_2| = r.
\label{eq:3.3}
\end{equation}
Let
\begin{equation}
\bar m := K - m \ \text{ for } m \in [0, K], \quad
n^* := \max(\sqrt K, n),\ \ \overline n^* := (\overline n)^* .
\label{eq:3.4}
\end{equation}
Our first main result is the following theorem.

\begin{theorem}
\label{ujth:1}
We have for $N^c < R \leq N^{1/2} \exp( - c \sqrt{\log N})$, as $N \to
\infty$
\begin{align}
&\sum_{\substack{n \leq N\\ n \in A(\mathcal H_1) \cap A(\mathcal
H_2)}} \Lambda_R(n;\ \mathcal H_1, \ell_1) \Lambda_R(n; \mathcal
H_2, \ell_2) =
\label{eq:3.5}\\
&= N {\ell_1 + \ell_2 \choose \ell_1} \frac{(\log R)^{r + \ell_1 +
\ell_2}}{(r + \ell_1 + \ell_2)!}
\frac{\mathfrak S(\mathcal H)P}
{|A(\mathcal H)|} \left(1 + O
\left( \frac{K \bar r^* \log_2 N}{\log R} \right)\right) +
O\left(N e^{-c\sqrt{\log N}}\right). \nonumber
\end{align}
\end{theorem}

For the next theorem we suppose that the following form of
the Bombieri--Vinogradov theorem holds (see
Section~\ref{sec:11}). For a given, sufficiently large $N$,
and recalling the parameter $P$ defined in \eqref{eq:2.14}, we have
\begin{equation}
\sum_{\substack{q \leq Q^*\\ (q, P) = 1}}
\max_{\substack{a\\(a,q) = 1}} \bigg| \sum_{\substack{N < p \leq
2N\\ p \equiv a(\text{\rm mod}\, Pq)}} \log p - \frac{N}{\varphi(Pq)}
\bigg| \ll \frac{N}{P} \exp (- c \sqrt{\log N}),
\label{eq:3.6}
\end{equation}
where
\begin{equation}
Q^* = N^{1/2} P^{-3} \exp \bigl(-c^* \sqrt{\log N}\bigr),
\label{eq:3.7}
\end{equation}
with an arbitrary positive constant $c^*$.

Letting $\mathcal H^0 = \mathcal H \cup \{h_0\}$, our second main result is as follows.

\begin{theorem}
\label{ujth:2}
Suppose \eqref{eq:3.6}--\eqref{eq:3.7} hold and let
$ N^c \leq R \leq \sqrt{Q^*}$.
Then
\begin{align}
&\sum_{\substack{N < n \leq 2N\\ n \in A(\mathcal H_1) \cap
A(\mathcal H_2)}}\Lambda_R(n;\mathcal{H}_1,\ell_1)
\Lambda_R(n;\mathcal{H}_2, \ell_2)\theta(n+h_0)
\label{eq:3.8}\\
&= N \frac{C_R(\ell_1, \ell_2, \mathcal{H}_1,\mathcal{H}_2,h_0) }
{(r + \ell_1 + \ell_2)!}
{\ell_1 + \ell_2 \choose \ell_1}
\gs({\mathcal{H}}^0)  (\log R)^{r + \ell_1 + \ell_2}
\left( 1 + O \left(\frac{CK \bar r^* \log_2 N}{\log R} \right)
\right)\nonumber\\
&\qquad  + O\left(N e^{- c\sqrt{\log N}}\right),
\nonumber
\end{align}
where
\begin{equation}
C_R(\ell_1, \ell_2, \mathcal{H}_1,\mathcal{H}_2,h_0) = \left\{
\begin{array}{ll}
      {1,} &
        \mbox{if $h_0\not \in \mathcal{H}$;}\\
 \frac{(\ell_1 + \ell_2 + 1) \log R}{(\ell_1 + 1)(r + \ell_1 + \ell_2 +1)},
 &
        \mbox{if $h_0 \in \mathcal{H}_1$ and $h_0\not \in
 \mathcal{H}_2$;}\\
 \frac{(\ell_1 + \ell_2 + 2)(\ell_1 + \ell_2 + 1) \log R}
{(\ell_1 + 1)(\ell_2 + 1)(r + \ell_1 + \ell_2 + 1)}, &
        \mbox{if $h_0 \in \mathcal{H}_1\cap \mathcal{H}_2$}.\\
\end{array}
\right.\label{eq:3.9}
\end{equation}
\end{theorem}

For the applications to
Theorems~\ref{th:A}--\ref{th:uj3} the simpler case $\ell_1 =
\ell_2 = \ell$ will be sufficient.

\section{Lemmas}
\label{sec:4}

We will use standard properties of the Riemann zeta function $\zeta(s)$.
Proceeding slightly differently from \cite{GPY} we use
the zero-free region, with $s = \sigma + it$,
\begin{equation}
\zeta(1 + s) \neq 0\  \text{ for } \ s \in \mathcal R_N := \left\{
s;\ \sigma \geq  - \frac{1}{\log_2 N + 6\log (|t| + 3)}\right\}.
\label{eq:4.1}
\end{equation}
Further we have for $s \in \mathcal R_N$
by Titchmarsh \cite[Ch.~3]{Tit}
\begin{equation}
\max \left(\left|\zeta(1 + s) - \frac{1}{s} \right|, \left|
\frac{1}{\zeta(1 + s)} \right|, \left| \frac{\zeta'}{\zeta}
(1 + s) +
\frac{1}{s} \right|\right) \ll \log (|t| + 3).
\label{eq:4.2}
\end{equation}

In the course of the proof the following contours 
which lie in the zero-free region $\mathcal R_N$
will be used (with $U$ and $V$ given in \eqref{eq:2.13})
\begin{equation}
\alignedat2
\mathcal L_1 &:= \left\{ {\sigma} = \frac1{28V}, \ |t| \leq
U  \right\}, &\quad
\mathcal L_2 &:= \left\{ {\sigma} = \frac1{14V}, \ |t| \leq
2U  \right\},\\
\mathcal L_3 &:= \left\{ {\sigma} = \frac{-1}{28V}, \ |t| \leq
U  \right\}, &\quad
\mathcal L_4 &:= \left\{ {\sigma} = \frac{-1}{14V}, \ |t| \leq
2U  \right\},\\
\mathcal L_5 &:= \left\{ - \frac{1}{28V} \leq {\sigma} \leq
\frac1{28V},  \ |t| = U  \right\}, &\quad
\mathcal L_6 &:= \left\{ - \frac1{14V} \leq {\sigma} \leq
\frac{1}{14V},\  |t| = 2U  \right\},
\endalignedat
\label{eq:4.3}
\end{equation}
\begin{align}
\mathcal L' & = \mathcal L'_0 \cup \mathcal L'_1; \quad
\mathcal L'_0 = \left\{ s; \ s = \delta_0 e^{i\varphi}, \
\frac{\pi}{2} \leq \varphi \leq \frac{3\pi}{2} \right\},
\label{eq:4.4}\\
\mathcal L'_1 &= \{ s = it; \ \delta_0 \leq |t| \leq U\},
\quad
\delta_0 = \left(\sqrt{K} \log_2 N\right)^{-1}. \nonumber
\end{align}

Similarly to Lemma~1 of \cite{GPY}, we have

\begin{lemma}
\label{l:1}
Let $k(\log_2 N)^2 \leq \sqrt{\log R}$, $N^c \leq R \leq N$,
$N \geq C$, $k \geq 2$, $B
\leq Ck$.
Then
\begin{equation}
\int\limits_{\mathcal L_i}(\log (|t| + 3))^B \left|\frac{R^s
ds}{s^k}\right| \ll e^{-c\sqrt{\log N}}, \quad
(3 \leq i \leq 6),
\label{eq:4.5}
\end{equation}
where the constant implied by the $\ll $ symbol depends only
on~$C$.
\end{lemma}

\noindent
{\it Proof.}
The integral $I$ in \eqref{eq:4.5} satisfies
for all $i$
\begin{align}
I &\ll \int\limits^{2U}_0 R^{- 1/(28V)}  \frac{(\log (|t| +
3))^B}{\max \left(|t|, \frac1{28V}\right)^k}\, dt
+ \int\limits_{|\sigma| \leq 1/14V} R^\sigma \cdot
\frac{d\sigma}{U^{3/2}} 
\label{eq:4.6}\\
&\ll e^{-c\sqrt{\log N}} \Bigg(\int\limits^C_0 (28V)^k \, dt +
\int\limits^\infty_C \frac{dt}{t^{3/2}} \Bigg) + e^{\sqrt{\log
N}(1/14 - 1/2)}  \ll e^{-c\sqrt{\log N}}.
\nonumber
\end{align}

We will prove a generalization of
the combinatorial identity (8.16)
of \cite{GPY} in order to evaluate the terms $I_{1,1}$ of
Section~\ref{sec:B}.
Let us define for triplets of integers $d, u, y$ with $d \geq
0$, $u \geq 0$, $y + u \geq 0$ (to be called suitable triplets) the
quantity
\begin{equation}
Z(d, u , y) := \frac{1}{u!} \sum^u_{\substack{m = 0\\ m \geq - y}} {u
\choose m} (-1)^m \frac{d(d + 1) \cdots (d + m - 1)}{(y + m)!}\, .
\label{eq:4.7}
\end{equation}

\begin{lemma}
\label{l:2}
We have for any suitable triplet $d, u, y$ the relation
\begin{equation}
Z(d, u, y) = \frac{(y - d + 1) \cdots (y - d + u)}{u!(y + u)!}\, .
\label{eq:4.8}
\end{equation}
\end{lemma}

\noindent
{\it Proof.}
We will prove this by induction on~$u$.
For $u = 0$ we have trivially for any non-negative $d$ and $y$,
$Z(d,0,y) = (y!)^{-1}$
(the empty product in the numerator of $Z$ is $1$ by definition).
We can suppose $u \geq 1$ and that our statement is true for
all suitable triplets $d$, $u - 1$, $y$.
Making the convention that we define for $n < 0$
\begin{equation}
\frac{x}{n!} = 0
\label{eq:4.9}
\end{equation}
for any real number $x$ (in other words, we just neglect in a
sum all terms with an $n!$ in the denominator with $n < 0$), we obtain
by ${u\choose i} = {u - 1\choose i} + {u - 1\choose i - 1}$
(where we define ${u - 1\choose u} = {u - 1\choose -1} = 0$),
with the notation $[S] = 1$ if the statement $S$ is true and
$[S] = 0$ if $S$ is false,
\begin{align}
Z(d,u,y)  \label{eq:50}
&= \frac{1}{u!} \bigg\{ \sum^{u - 1}_{\substack{i = 0\\ i \geq
- y}} {u -
1\choose i} (-1)^i \frac{d(d + 1) \dots (d + i - 1)}{(y + i)!}
\nonumber\\
&\quad - \sum^{u - 1}_{\substack{j = 0 \\ j \geq - y - 1}} {u - 1\choose
j} (-1)^j \frac{d(d + 1) \dots (d + j)}{(y + j + 1)!} \bigg\}
\nonumber\\
&= \frac{1}{u!} \bigg\{ \sum^{u - 1}_{\substack{i = 0\\ i \geq - y}} {u -
1\choose i} (-1)^i \frac{d(d + 1) \dots (d + i - 1)}{(y + i)!}
\left(1 - \frac{d + i}{y + i + 1}\right) \nonumber\\
&\quad \! -\!  [\! -\! y \! -\!  1 \geq 0] {u \! -\!  1\choose
\! -\!  y \! -\!  1} (\! -\! 1)^{-y - 1} d
(d\! +\! 1) \dots (d \! -\!  y \! -\!  2)(d \! -\!  y \! -\!  1) \bigg\}\nonumber \\
&= \frac1u \cdot \frac1{(u - 1)!}
\sum^{u - 1}_{\substack{i = 0\\ i
\geq - y - 1}} {u - 1\choose i} (-1)^i \frac{d(d + 1) \dots (d +
i - 1) (y + 1 - d)}{(y + i +1)!} \nonumber\\
&= \frac{y + 1 - d}{u} Z(d,u - 1, y + 1) = \frac{(y + 1 - d)(y +
2 - d) \dots (y + u - d)}{u!(y + u)!} .
\nonumber
\end{align}

Finally we mention a simple lemma for the mean value of the
generalized divisor function
\begin{equation}
d_m(q) := m^{\omega(q)},
\label{eq:uj4.11}
\end{equation}
where $\omega(q)$ denotes the number of prime-factors of $q$ for
a squarefree~$q$.

\begin{lemma}
\label{lem:uj2a}
If $m > 0$, $\nu \geq \max(c' \log (K + 1),
1)$ then there exists a constant $C'$ depending on $c'$ such
that, for $K \geq 1$ and $x \geq 1$ we have
\begin{equation}
\underset{q \leq x}{\sum\nolimits^\flat}\,  d_m(q) \leq x(1 +
\log x)^{\lceil m\rceil}
\label{eq:uj4.12}
\end{equation}
and
\begin{equation}
\underset{q \leq x}{\sum\nolimits^\flat}\, \frac{(d_{3K}(q))^{1
+ 1/\nu}}{q} \leq (1 + \log x)^{C'K}.
\label{eq:uj4.13}
\end{equation}
\end{lemma}

\noindent
{\it Proof.} Equation
\eqref{eq:uj4.12} follows from
\begin{equation}
\underset{q\leq x}{\sum\nolimits^\flat} d_m(q) \leq x \bigg(
\underset{q\leq x}{\sum\nolimits^\flat} \frac{d_{\lceil m
\rceil}(q)}{q} \bigg)
\leq x \bigg( \sum_{j \leq x} \frac1j\bigg)^{\lceil m \rceil}
\leq x(1 + \log x)^{\lceil m \rceil} .
\label{eq:uj5.14}
\end{equation}

Further, by \eqref{eq:uj4.11} we have
\begin{equation}
(d_{3K}(q))^{1+1/\nu} = d_j(q)
\label{eq:uj5.15}
\end{equation}
with
\begin{equation}
j = (3K)^{1 + 1/\nu} \leq 9 e^{1/c'} K,
\label{eq:uj5.16}
\end{equation}
and Lemma~\ref{lem:uj2a} follows with $C' = 9e^{1/c'} + 1$.

\section{Preparation for the Proof of Theorem~\ref{ujth:1} }
\label{sec:5}

Since the preparation for the proof of Theorem~\ref{ujth:1} is nearly
 the same as in Sections~6 and 7 of \cite{GPY} for the analogous
Proposition~1 (or 3 or 4), we will briefly summarize it
and the reader is referred for the details to \cite{GPY}.
Let
\begin{equation}
\mathcal H(p) = \left\{ h'_1, \dots , h'_{\nu_p(\mathcal H)} :
h'_j \equiv h_i\,(\text{\rm mod}\,p), \ h_i \in \mathcal H
\text{ for some } i, 1 \leq h'_j
\leq p\right\},
\label{eq:5.1}
\end{equation}

\begin{equation}
\bar \nu_p(\mathcal H_1 \bar \cap \mathcal H_2) := \nu_p
(\mathcal H_1(p) \cap \mathcal H_2(p)) = \nu_p (\mathcal H_1 ) +
\nu_p(\mathcal H_2) - \nu_p(\mathcal H).
\label{eq:5.2}
\end{equation}
For any $a \in A(\mathcal H_1) \cap A(\mathcal H_2)$
(cf.\ \eqref{eq:2.16}), we have similarly to Section 7 of \cite{GPY}
\begin{align}
S_R(N; \mathcal H_1, \mathcal H_2,
\ell_1, \ell_2, a)
:&=
\sum^{2N}_{\substack{n = N + 1\\
n\equiv a(\text{\rm mod}\, P)}} \Lambda_R(n; \mathcal H_1,
\ell_1) \Lambda_R(n; \mathcal H_2, \ell_2)
\label{eq:5.3}\\
&= \frac{N}{P}
\mathcal T_R(\ell_1, \ell_2; \mathcal H_1, \mathcal H_2) +
O(R^2(3 \log R)^{7K}),
\nonumber
\end{align}
where
\begin{equation}
\mathcal T_R(\ell_1, \ell_2; \mathcal H_1, \mathcal H_2) := \frac{1}{(2\pi
i)^2} \int\limits_{(1)}\!\int\limits_{(1)} F(s_1, s_2)
\frac{R^{s_1}}{s_1^{K + \ell_1 + 1}} \frac{R^{s_2}}{s^{K + \ell_2 +
1}_2}\, ds_1 ds_2,
\label{eq:5.4}
\end{equation}
\begin{equation}
F(s_1, s_2) := \prod_{p > V} \left( 1 - \frac{\nu_p (\mathcal
H_1)}{p^{1 + s_1}} - \frac{\nu_p(\mathcal H_2)}{p^{1 + s_2}} +
\frac{\bar \nu_p(\mathcal H_1 \bar \cap \mathcal H_2)}{p^{1 +
s_1 + s_2}} \right),
\label{eq:5.5}
\end{equation}
where now, differently from \cite{GPY}, primes not exceeding $V$ do
not appear in $F(s_1, s_2)$ since by  the regularity of $a$, $(P_{\mathcal H_1}(n), P) = (P_{\mathcal H_2}(n), P) =
1$.

Let
\begin{equation}
\Delta : = \Big|\prod_{1 \leq i < j \leq K} (h_i - h_j) \Big|
\leq N^{K(K - 1)/2}.
\label{eq:5.6}
\end{equation}
Then if $p \nmid \Delta$ (consequently, for all sufficiently
large primes $p$),
\begin{equation}
\nu_p(\mathcal H_1) = |\mathcal H_1| = K, \quad
\nu_p(\mathcal H_2) = |\mathcal H_2| = K, \quad
\bar \nu_p(\mathcal H_1 \bar \cap \mathcal H_2) = |\mathcal H_1
\cap \mathcal H_2| = r.
\label{eq:5.7}
\end{equation}
We therefore factor out the dominant zeta-factors and write
\begin{equation}
F(s_1, s_2) = G_{\mathcal H_1, \mathcal H_2}(s_1, s_2)
\frac{\zeta(1 + s_1 + s_2)^d}{\zeta(1 + s_1)^a \zeta(1 + s_2)^b}
\label{eq:5.8}
\end{equation}
with a function $G(s_1, s_2)$,
regular for $\sigma_i > -1/5$, say,
which we write slightly more
generally for future application in Theorem~\ref{ujth:2} as
\begin{equation}
G_{\mathcal H_1, \mathcal H_2}(s_1, s_2) = G(s_1, s_2) = G = G_1
G_2 G_3 = G_1 G_4,
\label{eq:5.9}
\end{equation}
where now $a = b = K$, $d = r$,
$\nu_1(p) = \nu_p(\mathcal
H_1)$, $\nu_2(p) = \nu_p(\mathcal H_2)$, $\nu_3(p) =
\bar\nu_p(\mathcal H_1 \bar\cap  \mathcal H_2)$,
\begin{equation}
G_1(s_1, s_2) = \prod_{p \leq V} \left(1 - \frac{1}{p^{1 + s_1}}
\right)^{-a} \prod_{p \leq V} \left(1 - \frac{1}{p^{1 +
s_2}}\right)^{-b} \prod_{p \leq V} \left(1 - \frac{1}{p^{1 + s_1
+ s_2}}\right)^d,
\label{eq:60}
\end{equation}
\begin{align}
G_4(s_1, s_2) &= \prod_{p > V} \left(\frac{\left(1 -
\frac{\nu_1(p)}{p^{1 + s_1}} - \frac{\nu_2(p)}{p^{1 + s_2}} +
\frac{\nu_3(p)}{p^{1 + s_1 + s_2}} \right) \left(1 -
\frac{1}{p^{1 + s_1 + s_2}}\right)^d}{\left(1 - \frac{1}{p^{1 +
s_1}} \right)^a \left(1 - \frac{1}{p^{1 + s_2}}\right)^b} \right)
\label{eq:61}\\
&= \prod_{p \mid \Delta, p > V} \cdot \prod_{p\nmid \Delta, p >
V} = : G_2(s_1, s_2) G_3(s_1, s_2).\nonumber
\end{align}

Let us use the notation
\begin{equation}
\delta_i := \max(0, - \sigma_i), \quad \delta := \delta_1 +
\delta_2, \quad s_3 := s_1 + s_2,
\label{eq:62}
\end{equation}
and
\begin{equation}
\mathcal R'_N := \left\{s;\ \sigma \geq -
\frac{1/2}{\log_2 N + 6\log(|t| + 3)}\right\}.
\label{eq:63}
\end{equation}
We will estimate the order of $G(s_1, s_2)$ in the region $s_1,
s_2 \in \mathcal R'_N$ under the more general conditions
\begin{equation}
a,b,d \leq K, \quad \nu_i (p) \leq K,
\label{eq:64}
\end{equation}
\begin{equation}
\nu_1 (p) = a, \quad \nu_2(p) = b, \quad \nu_3(p) = d \ \text{ for
} p \nmid \Delta.
\label{eq:65}
\end{equation}
(Later we will examine more delicate properties of $G(s_1, s_2)$
with further conditions on $a,b, d, \nu_i(p)$.)
We have
\begin{equation}
|G_1(s_1, s_2)| \leq \exp \bigg( C \sum_{p \leq V} \frac{K}{p^{1
- \delta}}\bigg) \leq \exp(CK \log_3 N),
\label{eq:66}
\end{equation}
\begin{align}
|G_2(s_1, s_2)| & \leq \exp \bigg( C \sum_{p \mid \Delta}
\frac{K}{p^{1 - \delta}}\bigg) \leq \exp \bigg(CK \sum_{p \leq
\log \Delta (1 + o(1))} \frac1{p^{1 - \delta}} \bigg)
\label{eq:67}\\
&\leq \exp (CK\log_3 N), \nonumber
\end{align}
and
\begin{equation}
|G_3(s_1, s_2)| \leq \exp \bigg(C \sum_{p > V} \frac{K^2}{p^{2 -
2\delta}}\bigg) \leq \exp \left(\frac{CK^2}{V}\right) \leq \exp (CK),
\label{eq:68}
\end{equation}
where in \eqref{eq:66}--\eqref{eq:68} we made use of the estimates
\begin{equation}
\max\big(V^\delta, (\log\Delta)^\delta\big) \leq (\log^2
N)^{\frac{1/2}{\log_2 N}} = e;
\label{eq:69}
\end{equation}
further in \eqref{eq:67} the sum which was originally over $p \mid
\Delta$ has been majorized by using the set of the smallest
possible primes which could divide~$\Delta$.

Summarizing \eqref{eq:66}--\eqref{eq:68} we obtain
\begin{align}
|G(s_1, s_2)| &\leq e^{CK \log_3 N} \quad \text{ for } \
s_1, s_2 \in \mathcal R'_n,
\label{eq:71}\\
\intertext{and further}
|F(s_1, s_2)| &\leq e^{CK \log_3 N} ((\log(|t_1| + 3)) \log (|t_2|
+ 3))^{2K} \text{ for } s_1, s_2, s_3 \in \mathcal R'_n.
\label{eq:72}
\end{align}

The above estimate shows that the integrand  in
\eqref{eq:5.4} vanishes as either $|t_1| \to \infty$ or $|t_2| \to
\infty$, $s_1, s_2, s_3 \in {\mathcal R}'_N$.
We will examine the integral (analogously to (7.15) in \cite{GPY})
\begin{equation}
I:= \mathcal T^*_R(d,a,b,u,v,\mathcal H_1, \mathcal H_2) :=
\frac1{(2\pi i)^2} \int\limits_{(1)}\! \int\limits_{(1)}
\frac{D(s_1, s_2) R^{s_1 + s_2} ds_1 ds_2}{s^{u + 1}_1 s^{v +
1}_2 (s_1 + s_2)^d},
\label{eq:73}
\end{equation}
where we introduce the function $D(s_1, s_2)$, regular for $s_1,
s_2, s_3 \in {\mathcal R}'_N$,
\begin{equation}
D(s_1, s_2) := D_0(s_1, s_2)G(s_1, s_2), \
D_0(s_1, s_2) := \frac{W^d(s_1 + s_2)}{W^a(s_1) W^b(s_2)}, \
W(s)
:= s\zeta(1 + s),
\label{eq:74}
\end{equation}
\begin{equation}
0 \leq d \leq a, \ b\leq K,\ \ \min(a,b) \geq cK, \ \ \sqrt K /
8 \leq u, \ v \leq \sqrt K
\label{eq:75}
\end{equation}
and by symmetry we can assume $u \leq v$. (In the applications
we will have $|a - b| \leq 1$, $|u - v| \leq 1$.)

{\it First step.}
Move the contour $(1)$ for the integral over $s_1$ to $\mathcal
L_1$,
over $s_2$ to $\mathcal L_2$.
The vertical parts $|t| \geq U$, and $|t| \geq 2U$, resp.\ can
be neglected similarly to Lemma~1.
After this move the integral in $s_1$ from $\mathcal L_1$ to
$\mathcal L_3 \cup \mathcal L_5$.
The horizontal segments $\mathcal L_5$ can be again neglected.
We pass a pole of order $u + 1$ at $s_1 = 0$, and obtain
\begin{equation}
I = I_1 + \frac{1}{(2\pi i)^2} \int_{\mathcal L_2} \int_{\mathcal L_3}
\frac{D(s_1, s_2)
R^{s_1 + s_2} ds_1 ds_2}{s^{u+1}_1 s^{v+1}_2 (s_1+s_2)^d} = I_1
+ I_2 + O\Bigl(e^{-c\sqrt{\log N}}\Bigr),
\label{eq:5.18}
\end{equation}
where
\begin{equation}
\aligned
I_1 :&= \frac{1}{2\pi i} \int_{\mathcal L_2} \text{\rm Res}_{s_1 = 0}
\left(\frac{D(s_1, s_2) R^{s_1 + s_2}}{s^{u+1}_1 s^{v+1}_2
(s_1+s_2)^d} \right) ds_2\\
&= \frac{1}{2\pi i}
\int\limits_{\mathcal L_2} \frac1{u!} \bigg\{
\sum^u_{i=0} {u\choose i} (\log R)^{u-i}
\frac{\partial^i}{\partial s^i_1} \left( \frac{D(s_1, s_2)}{(s_1
+ s_2)^d}\right) \bigg|_{s_1 = 0} \bigg\}
\frac{R^{s_2}}{s^{v+1}_2} ds_2.
\endaligned
\label{eq:5.19}
\end{equation}
We denote the complete integrand above by $Z(s_2)$ and express
\begin{equation}
\gathered
\frac{\partial^i}{\partial s^i_1} \left(\frac{D(s_1, s_2)}{(s_1
+ s_2)^d}\right) \bigg|_{s_1 = 0} = (-1)^i \frac{D(0, s_2)
d(d+1) \dots (d + i - 1)}{s^{d+i}_2}
\\
+ \sum^i_{j = 1} {i\choose j} \frac{\partial^j}{\partial s^j_1} D(s_1, s_2)
\bigg|_{s_1 = 0} \cdot (-1)^{i-j} \frac{d(d+1) \dots (d + i - j -
1)}{s^{d+i - j}_2}
\endgathered
\label{eq:5.20}
\end{equation}
where in case of $i = j$
(including also the case when $i = j = 0$ and $d \geq  0$ arbitrary)
the empty product  in the numerator is $1$.

{\it Second step.}
Let us denote the contribution of the first term in \eqref{eq:5.20} to
\eqref{eq:5.19} by $I_1(i, 0)$ and the others by $I_1(i,j)$ $(1 \leq j
\leq i)$. $I_1(i,0)$ will belong to the main term,
all $I_1(i,j)$ with $j \geq 1$ will just contribute to the secondary terms.
Let us move now the contour $\mathcal L_2$ for the integral over $s_2$
to $\mathcal L_4 \cup \mathcal L_6$ in~\eqref{eq:5.19}.
The horizontal segments $\mathcal L_6$ can be neglected again.
We pass a pole of order $v + 1 + d + i - j$ in case of
$I_1(i,j)$ and we obtain in this way
\begin{equation}
\gathered
I_1 = \frac1{u!} \sum^u_{i=0} {u\choose i} (\log R)^{u - i}
\sum^i_{j=0} (-1)^{i-j}
{i\choose j} \frac{d(d+1)\dots (d+i - j -
1)}{(v+d+i-j)!} \times\\
\times \sum^{v+d+i-j}_{\nu=0} {v+d+i-j \choose \nu} (\log
R)^{v+d+i-j-\nu} \cdot \frac{\partial^\nu }{\partial s^\nu_2}
\frac{\partial^j}{\partial s^j_1} D(s_1, s_2) \bigg|_{s_1 = s_2
= 0} \\
+ \frac{1}{2\pi i} \int\limits_{\mathcal L_4} Z(s_2) ds_2
+ O\bigl(e^{-c\sqrt{\log N}}\bigr)
= :
I_{1,1} + I_{1,2} + O\bigl(e^{-c\sqrt{\log N}}\bigr).
\endgathered
\label{eq:6.23}
\end{equation}

\section{Estimates of the partial derivatives of $D(s_1,s_2)$}
\label{sec:A}

In this section we will estimate partial derivatives
$\frac{\partial^i}{\partial s^i_1} \frac{\partial^j}{\partial
s^j_2} D(s_1, s_2)$ of $D(s_1, s_2)$ for $i + j \leq CK$ with
$s_i = s^*_i$ in
${\mathcal R}'_N$ for $1 \leq i \leq 3$.
We will often use Cauchy's estimate for functions regular in $|z
- z_0| \leq \eta$:
\begin{equation}
\frac1{j!} |f^{(j)} (z_0)| \leq \eta^{-j} \max_{|z-z_0| = \eta} |f(z)|.
\label{eq:A1}
\end{equation}

Applying this for $D(s_1, s_2)$ we obtain
\begin{equation}
\frac1{i!j!} \left|\frac{\partial^i}{\partial s^i_1}
\frac{\partial^j}{\partial s^j_2} D(s^*_1, s^*_2)\right| \ll
\eta^{-(i + j)} \max_{|s'_1 - s^*_1| \leq \eta, |s'_2 - s^*_2|
\leq \eta} |D(s'_1, s'_2)|.
\label{eq:A2}
\end{equation}

In order to substitute the above maximum for $D(s^*_1,
s^*_2)$, we have to estimate
\begin{equation}
L(s_1, s_2) : = \max \left(\left| \frac{\partial}{\partial s_1}
\log D(s_1, s_2) \right|, \left|\frac{\partial}{\partial s_2}
\log D(s_1, s_2)\right| \right),
\label{eq:A3}
\end{equation}
for $s_1, s_2, s_3 \in {\mathcal R}'_N$; since by the regularity of $\log D(s_1,
s_2)$
for $s_i \in {\mathcal R}'_N$ $(1 \leq i \leq3)$
(cf.\ \eqref{eq:60}, \eqref{eq:61},
\eqref{eq:74}) we have, for $\eta \leq (\log_2 N +
\log(|t_1| + 3) + \log(|t_2| + 3))^{-1}/100$,
\begin{equation}
\left| \frac{D(s'_1, s'_2)}{D(s^*_1, s^*_2)}\right| \leq \exp
\big( 2\eta \cdot \max_{|s_1 - s^*_1| \leq \eta, |s_2 - s^*_2|
\leq \eta} L(s_1, s_2) \big).
\label{eq:A4}
\end{equation}

By symmetry it is enough to deal with
\begin{equation}
L_1(s_1, s_2) := \left| \frac{\partial}{\partial s_1} \log
D(s_1, s_2)\right|.
\label{eq:A5}
\end{equation}
Since the logarithm is an additive function, using the
representation \eqref{eq:5.9}--\eqref{eq:61} and \eqref{eq:74} of
$D(s_1, s_2)$ it is sufficient to examine the factors $D_0$,
$G_1$, $G_2$, $G_3$ separately.

We will choose a positive $\eta$
\begin{equation}
\eta \leq \frac{1/100}{\log_2 N + \log T}, \quad
T = T_1 + T_2 ,
\quad T_i = |t_i| + 3, \quad
t_3 = t_1 + t_2 \quad (1 \leq i \leq 3)
\label{eq:A6}
\end{equation}
where by $\delta = \delta_1 + \delta_2 \leq 2 / \log_2 N$, $V =
(\log N)^{1/2}$, $\log \Delta \leq K^2 \log N \leq \log^2 N$ we have
\begin{equation}
\max\big(V^\delta, V^\eta, (\log \Delta)^\delta, (\log
\Delta)^\eta\big) \ll 1.
\label{eq:A7}
\end{equation}

We have by \eqref{eq:4.2} and \eqref{eq:74}
\begin{align}
\frac{\partial}{\partial s_1} (\log D_0(s_1, s_2))
&= d \cdot
\left(\frac{\zeta'}{\zeta} (1 + s_1 + s_2) + \frac1{s_1 +
s_2}\right) - a \left( \frac{\zeta'}{\zeta} (1 + s_1) +
\frac1{s_1} \right)
\label{eq:A8}\\
&\ll K \log T.\nonumber
\end{align}
Further we have
\begin{equation}
\frac{\partial}{\partial s_1} (\log G_1(s_1, s_2)) =
\sum_{p \leq V} \frac{\log p}{p^{1 + s_1}}
\left(\frac{dp^{-s_2}}{1 - p^{-(1 + s_1 + s_2)}} - \frac{a}{1 -
p^{-(1 + s_1)}} \right) \ll K \log_2 N.
\label{eq:A9}
\end{equation}

Similarly to \eqref{eq:67} we obtain by \eqref{eq:A7}
\begin{equation}\begin{split}
&\frac{\partial}{\partial s_1}(\log G_2(s_1, s_2))
\\ &=
\sum_{\substack{p\mid \Delta\\ p > V}} \frac{\log p}{p^{1 + s_1}}
\bigg\{ \frac{\nu_1(p) - p^{-s_2} \nu_3(p)}{1 -
\frac{\nu_1(p)}{p^{1 + s_1}} - \frac{\nu_2(p)}{p^{1 + s_2}} +
\frac{\nu_3(p)}{p^{1 + s_1 + s_2}}} -
\frac{a}{1 - \frac{1}{p^{1 + s_1}}} +
\frac{dp^{-s_2}}{1 - \frac{1}{p^{1 + s_1 + s_2}}} \bigg\}\\
&\ll K \sum_{p \mid \Delta} \frac{\log p}{p^{1 - \delta}} \ll K
\sum_{p \leq \log \Delta(1 + o(1))} \frac{\log p}{p^{1 -
\delta}} \ll K \log_2 \Delta \ll K \log_2 N. 
\label{eq:A10}\end{split}\end{equation}

Finally, analogously to \eqref{eq:68} we have by \eqref{eq:65}
\begin{equation}\begin{split}
& \frac{\partial}{\partial s_1}(\log G_3(s_1, s_2))\\& =
\sum_{\substack{p\nmid \Delta\\ p > V}} \frac{\log p}{p^{1 + s_1}}
\bigg\{ \frac{ a - dp^{-s_2}}
{1 - \frac{a}{p^{1 + s_1}} - \frac{b}{p^{1 + s_2}} +
\frac{d}{p^{1 + s_1 + s_2}}} -
 \frac{a}{1 - \frac{1}{p^{1 + s_1}}}
+ \frac{dp^{-s_2}}{1 - \frac{1}{p^{1 +
s_1 + s_2}}} \bigg\} \\
& \ll \sum_{p > V} \frac{\log p}{p^{1 - \delta_1}}
\left(a \cdot \frac{K}{p^{1 - \delta}} + dp^{\delta_2} \cdot
\frac{K}{p^{1 - \delta}} \right)\ll K^2 \sum_{p > V} \frac{\log p}{p^{2 - 2\delta}} \ll
\frac{K^2}{V^{1 - 2\delta}} \ll \frac{K^2}{V} \ll K.
\label{eq:A11}\end{split}
\end{equation}

Summarizing \eqref{eq:A3}--\eqref{eq:A11} we have
\begin{equation}
\max_{|s'_1 - s^*_1| \leq \eta, |s'_2 - s^*_2| \leq \eta}
|D(s'_1, s'_2)| \leq e^{C\eta K(\log_2 N + \log T)} |D(s^*_1,
s^*_2)| \text{ if } s^*_1, s^*_2, s^*_3 \in \mathcal R'_N.
\label{eq:A12}
\end{equation}

Hence, \eqref{eq:A2} and \eqref{eq:A12} imply by the
choice $\eta^{-1} = 100 K(\log_2
N + \log T)$ the following estimate.

\begin{lemma}
\label{l:A1}
We have for $s_1, s_2, s_3 \in
\mathcal R'_N$
\begin{equation}
\frac1{i!j!} \left| \frac{\partial^i}{\partial s^i_1}
\frac{\partial^j}{\partial s^j_2} D(s_1, s_2) \right| \ll
\left(CK  (\log_2 N + \log T)\right)^{i + j} |D(s_1, s_2)|.
\label{eq:A13}
\end{equation}
\end{lemma}

The above estimate is sufficient for our purposes at every point
$(s_1, s_2)$ apart from $(0,0)$, which will appear in the main term.
We will show an analogous result for the point  $(0,0)$ where
$\eta$ in \eqref{eq:A6} will be replaced by the larger value
\begin{equation}
\eta_0 = \frac1{\phantom{\int^*}\!\!\!\bar d^* \log_2 N} .
\label{eq:A14}
\end{equation}
where we use the notation $\bar d, \bar d^*$ of \eqref{eq:3.4}.
Next we have
\begin{lemma}
\label{l:A2}
$\frac{1}{i!j!} \left| \frac{\partial}{\partial s^i_1}
\frac{\partial^j}{\partial s^j_2} D(s_1, s_2)\right|_{s_1 = s_2
= 0} \ll (\bar d^* \log_2 N)^{i + j} D(0,0)$.
\end{lemma}

\noindent
{\it Proof.}
Let $d_1 = a - d$.  Analogously to \eqref{eq:A8}--\eqref{eq:A9},
we have, for $|s_1|, |s_2| \leq \eta_0$, 
\begin{equation}\begin{split}
\frac{\partial}{\partial s_1} \log D_0(s_1, s_2) &= d
\frac{W'}{W} (s_1 + s_2) - a\frac{W'}{W}(s_1) \\&
= d \left(\frac{W'}{W}(s_1 + s_2) - \frac{W'}{W}(s_1)\right) -
(a - d) \frac{W'}{W}(s_1) \\
&\ll K\eta_0 +  d_1 \ll \frac{K}{\bar d^*} + d_1 \ll
\sqrt{K} +  d_1 \ll \bar d^*, \label{eq:A15}
\end{split}\end{equation}
and
\begin{equation}\begin{split}
\frac{\partial}{\partial s_1} (\log G_1(s_1&, s_2)) 
= \sum_{p \leq V} \frac{\log p}{p^{1 + s_1}}
\left(\frac{d}{p^{s_2} -  p^{-(1 + s_1)}} - \frac{d}{1 -
p^{-(1 + s_1)}} + \frac{d - a}{1 -  p^{-(1 + s_1)}} \right) \\
&\ll \sum_{p \leq V} \frac{\log p}{p^{1 - \delta_1}}
\left(\frac{K | p^{s_2} - 1|}{p^{-\delta_2}} + d_1 \right) \ll K \sum_{p \leq V} \frac{\eta_0 \log^2 p}{p^{1 - \delta}} +
 d_1 \sum_{p \leq V} \frac{\log p}{p^{1 - \delta}} \\
&\ll V^\delta \log V (K\eta_0 \log V + d_1) \ll
K\eta_0 \log^2_2 N + d_1 \log_2 N \ll \bar d^* \log_2 N.
\label{eq:A16}\end{split}
\end{equation}

The treatment of $G_2$ will be similar to this and \eqref{eq:A10}.
By $|\nu_1(p) - \nu_3(p)| = |\nu_p(\mathcal H_1) -
 \nu_p(\mathcal H_1(p) \cap \mathcal H_2(p))| \leq |\mathcal H_1
\setminus \mathcal H_2| = a - d$ we have
\begin{equation}\begin{split}
&\frac{\partial}{\partial s_1} (\log G_2(s_1, s_2)) =\\& \sum_{\substack{p
\mid \Delta\\ p > V}} \frac{\log p}{p^{1 + s_1}} 
\left\{ \frac{\nu_1(p) - \nu_3(p) + \nu_3(p)(1 -
 p^{-s_2})}{1 - \frac{\nu_1(p)}{p^{1 + s_1}} -
\frac{\nu_2(p)}{p^{1 + s_2}} + \frac{\nu_3(p)}{p^{1 + s_1 +
s_2}}} - \frac{d_1}{1 - \frac{1}{p^{1 + s_1}}}\right. \left. - d \left(
\frac1{1 - \frac1{p^{1 + s_1}}} - \frac{1}{p^{s_2} -
\frac{1}{p^{1 + s_1}}} \right)\right\} \\
&\ll \sum_{p\mid \Delta, p > V} \frac{\log p}{p^{1 - \delta}}
\left(d_1 + K(p^{\eta_0} - 1) \right)\ \ll  \ d_1 \sum_{p \mid \Delta} \frac{\log p}{p^{1 - \delta}}
+ K \eta_0 \sum_{\substack{p \mid \Delta \\ p \leq e^{2/\eta_0}}}
\frac{\log^2 p}{p^{1 - \delta}} + K \sum_{p \mid \Delta}
\frac{\log p}{e^{1/\eta_0}}\\
&\ll d_1 \sum_{p \leq \log \Delta(1 + o(1))} \frac{\log
p}{p^{1 - \delta}} + K \eta_0 \sum_{p \leq \log \Delta(1 +
o(1))} \frac{\log^2 p }{p^{1 - \delta}} + \frac{K \log
\Delta}{(\log N)^{\bar d^*}} \\
&\ll d_1 \log_2 \Delta + K \eta_0 \log^2_2 \Delta +
\frac1{(\log N)^{\sqrt K - 3}} \ll \log_2 N (d_1 + K \eta_0 \log_2 N) \ll \bar d^* \log_2
N.
\label{eq:A17} \end{split}\end{equation}
Finally we have, similarly to  above and \eqref{eq:A11},
using $a - dp^{-s_2} = a - d + d\bigl(1 - p^{-s_2}\bigr)$,
\begin{align}
&\frac{\partial}{\partial s_1} \log G_3(s_1, s_2)
= \sum_{\substack{p \nmid \Delta\\ p > V}} \frac{\log p}{p^{1 + s_1}}
\Bigg\{ a\left( \frac1{1 - \frac1{p^{1 + s_1 + s_2}}} - \frac1{1
- \frac1{p^{1 + s_1 }}} \right)
\label{eq:A18}\\
&\quad + (a - dp^{-s_2}) \left( \frac1{1 - \frac{a}{p^{1 + s_1}}
- \frac{b}{p^{1 + s_2}} + \frac{d}{p^{1 + s_1 + s_2}}} -
\frac1{1 - \frac1{p^{1 + s_1 + s_2}}}\right) \Bigg\} \nonumber\\
&\ll \sum_{p > V} \frac{\log p}{p^{1 - \delta_1}} \left\{
\frac{K}{p^{1 - \delta}} + \frac{K}{p^{1 - \delta}} \big( d_1 +
K(p^{\eta_0} - 1)\big)\right\}\nonumber\\
&\ll K d_1 \sum_{p > V}
\frac{\log p}{p^{2 - 2\delta}} + K^2\eta_0
\sum_{V < p \leq e^{1/\eta_0}} \frac{\log^2 p}{p^{2 - 2\delta}}
+ K^2 \sum_{p > e^{1/\eta_0}} \frac{\log p}{p^{2 - 2\delta -
\eta_0}} \nonumber\\
&\ll \frac{K d_1}{V^{1 - 2\delta}} + \frac{K^2 \eta_0 \log
V}{V^{1 - 2\delta}} + \frac{K^2}{(\log N)^{(1 -
2\delta - \eta_0)\bar d^*}} \nonumber\\
&\ll \frac{K}{V} ( d_1 + K\eta_0 \log_2 N) + o(1) \ll
d_1 + K \eta_0 \log_2 N \ll \bar d^*.
\nonumber
\end{align}

Now, \eqref{eq:A15}--\eqref{eq:A18} imply, by symmetry for $i = 1,2$, that
\begin{equation}
\left|\frac{\partial}{\partial s_i} \log D(s_1, s_2) \right| \ll
\eta^{-1}_0, \quad \text{for } |s_1|, |s_2| \leq \eta_0,
\label{eq:A19}
\end{equation}
and therefore, similarly to \eqref{eq:A12} we have
\begin{equation}
\max_{|s'_1| \leq \eta_0, |s'_2| \leq \eta_0} |D(s'_1, s'_2)|
\ll D(0,0),
\label{eq:A20}
\end{equation}
which by \eqref{eq:A2} proves Lemma~\ref{l:A2}.

\section{Contribution of the residue at $s_1 = s_2 = 0$}
\label{sec:B}

This section will be devoted to the examination of $I_{1,1}$,
the sum of the residues in \eqref{eq:6.23}.

The rather complicated formula \eqref{eq:6.23}
yields the main term and all secondary terms of the form $(\log
R)^m$
exclusively for $m \in [d, d+u+v-1]$ and will additionally
contribute to other secondary terms for $m \in [0, d - 1]$.
However, from the terms $I_{1,1}(i, j, \nu)$ belonging to
the triplet $(i,j,\nu)$ in the triple
summation, only those with $\nu = 0$, $j = 0$ contribute to the
main term of order $(\log R)^{d+u+v}$, since in all other terms
the exponent of $\log R$ is $d + u + v - j - \nu$.

We have to work now more carefully than in \cite{GPY}.
For example, by the aid of Lemma~\ref{l:2} (a generalization of
(8.16) in \cite{GPY}) we will exactly evaluate the
coefficients $A_{j,\nu}$ of $\frac1{j!\nu!}
\frac{\partial^\nu}{{\partial s_2}^\nu}
\frac{\partial^j}{{\partial s_1}^j} D
(s_1,
s_2)(\log R)^{v + d + u - j - \nu}$ in \eqref{eq:6.23} as follows.
Let $j, \nu \geq 0$,
\begin{equation}
m := i - j \geq 0, \quad  \  y : = v + d -
\nu,
\label{eq:B1}
\end{equation}
where we can assume by \eqref{eq:6.23}
\begin{equation}
\nu \leq v + d + m \Longleftrightarrow m \geq \nu - v - d = -y.
\label{eq:B2}
\end{equation}

Then we have from \eqref{eq:6.23},
by notation \eqref{eq:4.7}, \eqref{eq:B1} and Lemma~\ref{l:2}
\begin{equation}
A_{j,\nu} = \frac{j!\nu!}{u!} \sum^{u-j}_{\substack{m = 0\\ m \geq - y}}
{u \choose m + j} (-1)^m {m + j\choose j} \frac{d(d + 1) \dots
(d + m - 1)}{(v + d + m - \nu)! \nu!}
\label{eq:B3}
\end{equation}
\begin{align*}
&= \sum^{u - j}_{\substack{m = 0 \\ m \geq -y}}
\frac{(-1)^m}{(u - j - m)! m
!} \cdot \frac{d(d + 1) \dots (d + m - 1)}{(v + d + m - \nu)!}\\
&= Z(d, u - j, v + d - \nu)\ =\ \frac{(v - \nu + 1) \dots (v - \nu + u - j)}{(u - j)!(d + v -
\nu + u - j)!}\cdot
\end{align*}

We have to compare $A_{j,\nu}$ with $A_{0,0}$.
This will be furnished by the following

\begin{lemma}
\label{l:B1}
$|A'_{j,\nu}| := \left|\frac{A_{j,\nu}}{A_{0,0}}\right| \leq
(CK)^{j + \nu}$.
\end{lemma}

\noindent
{\it Proof.}
$|A'_{j,\nu}| = \frac{(d + v + u)!}{(d + v + u - \nu - j)!} \cdot
\frac{(u - j + 1) \dots u}{(v + u - j + 1) \dots (v + u)} \cdot
|A''_{j,\nu} | \leq (CK)^{j + \nu} |A''_{j,\nu}|$,
where
\begin{equation}
|A''_{j,\nu}| = \frac{|(v - \nu + 1) \dots (v - \nu + u -
j)|}{(v + 1) \dots (v + u - j)} .
\label{eq:B4}
\end{equation}

If $\nu \leq 2(v + 1)$, then clearly $A''_{j,\nu} \leq 1$, so we
may suppose
\begin{equation}
\nu = B(v + 1), \quad B > 2.
\end{equation}
In this case we have by $u - j \leq u \leq v < v + 1$:
\begin{equation}
A''_{j,\nu} \leq \left(\frac{\nu}{v + 1}\right)^{u - j} \leq
B^{v + 1} = B^{\nu/B} < 2^\nu,
\label{eq:B6}
\end{equation}
since the maximum of $x^{1/x}$ in $[1, \infty)$ is attained at
$x = e$ and $e^{1/e} < 2$.

\medskip
Now we are ready to evaluate the crucial term $I_{1,1}$ by the
aid of Lemmas~\ref{l:2}, \ref{l:A2} and \ref{l:B1}.
Namely by
\eqref{eq:3.1}, \eqref{eq:3.4}, $R \gg N^c$,
\eqref{eq:6.23}, \eqref{eq:B3} and \eqref{eq:75} we have
\begin{align}
I_{1,1} &= A_{0,0}(\log R)^{d + u + v}
\bigg\{ D(0,0) + \sum^u_{j
= 0} \sum^{v + d + u - j}_{\substack{\nu = 0\\ j + \nu \geq 1}}
\frac{A'_{j,\nu}}{(\log R)^{j + \nu}} \cdot
\frac{\frac{\partial^j}{\partial s_1}
\frac{\partial^\nu}{\partial s_2}
 D(0,0)}{j!\nu!}\bigg\}
\label{eq:B7}\\
&= Z(d, u, v + d) (\log R)^{d + u + v} D(0,0) \Bigg(\!1\! +
\! O \bigg(
\sum^\infty_{j = 0}  \sum^\infty_{\substack{\nu = 0 \\ j + \nu \geq 1}}
\left( \frac{CK\bar d^* \log_2 N}{\log R} \right)^{j + \nu}
\bigg)\! \Bigg)\nonumber \\
&= \frac{{v + u\choose u} (\log R)^{d + v + u}D(0,0)}{(d + v +
u)!} \left(1 + O \left(\frac{CK \bar d^* \log_2 N}{\log R}
\right) \right).\nonumber
\end{align}

The integral $I_{1,2}$ in \eqref{eq:6.23} does not contribute to the main
term and can be estimated relatively easily due to the presence
of the term $R^{s_2}$ $(s_2 \in \mathcal L_4)$.
In fact, choosing $\eta^{-1} = 100 (\log_2 N + \log T)$ as earlier, we
obtain by Lemma~\ref{l:A1}, \eqref{eq:4.2}, \eqref{eq:71},
\eqref{eq:5.19}--\eqref{eq:5.20}
for any $s_2\in \mathcal L_4$,
\begin{align}
Z(s_2) &\ll \sum^u_{i = 0} \sum^i_{j = 0} \frac{(\log R)^{u - i}
(CK)^{i - j}}{(u - i)!(i - j)!} (CK(\log_2 N + \log T))^j
\frac{|D(0, s_2)|R^{\sigma_2}}{|s_2|^{d + i - j + v + 1}}
\label{eq:B8}\\
&\ll e^{C \sqrt K \log_2 N} \frac{(\log(|t_2| + 3))^{3K + O(\sqrt{K})}
R^{\sigma_2}}{|s_2|^{b + O(\sqrt{K})}} .
\nonumber
\end{align}
Now Lemma~\ref{l:1} yields immediately by \eqref{eq:75}
and $R \gg N^c$
\begin{align}
I_{1,2}
&= \frac1{2\pi i} \int_{\mathcal L_4} Z(s_2) ds_2
\ll e^{C \sqrt{K} \log_2 N -
 c\sqrt{\log N}}
\label{eq:B9}\\
&\ll e^{- c \sqrt{\log N}}.
\nonumber
\end{align}

We may summarize \eqref{eq:B7} and \eqref{eq:B9} by
$D(0,0) = D_0(0,0) G(0,0) = G(0,0) \neq 0$
(which is true by the admissibility condition)
by
\eqref{eq:5.9}--\eqref{eq:61} and \eqref{eq:74} as

\begin{lemma}
\label{l:B2}
The integral $I_1$ in \eqref{eq:5.19} satisfies the asymptotic
\begin{equation}
I_1 = \frac{{v + u\choose u} (\log R)^{d + v + u} G(0,0)}{(d + v +
u)!} \left(1 + O \left( \frac{K \bar d^* \log_2 N}{\log R}
\right)\right) + O \big(e^{- c\sqrt{\log N}}\big).
\label{eq:B10}
\end{equation}
\end{lemma}

\section{Estimate of the integral $I_2$}
\label{sec:C}

For $I_2$ in \eqref{eq:5.18}, after interchange of the two integrations
 we move the contour $\mathcal L_2$  for the inner integral over $s_2$
to the left to
$\mathcal{L}_4$ passing a pole of order $d$ at $s_2= -s_1$
if $|t_2| \leq U$
and a
pole of order $v + 1$ at $s_2=0$ and obtain
\begin{equation}\begin{split}
 I_2 &= \frac {1}{ 2\pi
i}\mathop{\int}_{\mathcal{L}_3}
\mathop{\mathrm{Res}}_{s_2=-s_1}
\Big(\frac{D(s_1, s_2)
R^{s_1 + s_2} }{s^{u+1}_1 s^{v+1}_2 (s_1+s_2)^d}
\Big)\, ds_1
\label{eq:6.42} + \frac {1}{ 2\pi i}\mathop{\int}_{\mathcal{L}_3}
\mathop{\mathrm{Res}}_{s_2 = 0}
\Big(\frac{D(s_1, s_2)
R^{s_1 + s_2} }{s^{u+1}_1 s^{v+1}_2 (s_1+s_2)^d}
\Big)\, ds_1 \\& \qquad  + \frac {1}{ (2\pi
i)^{2}}\mathop{\int }_{\mathcal{L}_4}\mathop{\int}_{\mathcal{L}_3}
F(s_1,s_2)\frac{R^{s_1}}{ {s_1}^{a + u +1}} \frac{R^{s_2}}{
{s_2}^{b + v +1}}ds_1ds_2
+ O \bigl(e^{-c \sqrt{\log N}}\bigr) 
\\& := I_{2,1} + I_{2,2} + I_{2,3}
 + O \bigl(e^{-c \sqrt{\log N}}\bigr).
\end{split}\end{equation}
By the argument of
Lemma~\ref{l:1} and \eqref{eq:72},
the third integral $I_{2,3}$ is
$\ll e^{-c\sqrt{\log N}}$.
The second integral $I_{2,2}$ is completely analogous to
$I_{1,2}$ in \eqref{eq:6.23}, which was estimated by
$e^{-c\sqrt{\log N}}$ in \eqref{eq:B9}, the only change
being that the role of $s_1$ and $s_2$ is interchanged.

The residue in $I_{2,1}$ is zero if $d = 0$, while
for $d\ge 1$ we have 
\begin{equation}\begin{split}
\mathop{\mathrm{Res}}_{s_2=-s_1}\Big(
\frac{D(s_1, s_2)
R^{s_1 + s_2} }{s^{u+1}_1 s^{v+1}_2 (s_1+s_2)^d}
\Big)
&= \lim_{s_2\to
-s_1}\frac{1}{(d-1)!}\frac{\partial^{d-1}}{{\partial
s_2}^{d-1}}\left(\frac{
D(s_1,s_2)   R^{s_1+s_2}}{{s_1}^{u+1}
{s_2}^{v+1} }\right)\\
& = \frac1{(d-1)!}
\sum^{d-1}_{j=0} \mathcal{B}_j(s_1,\mathcal{H}_1,\mathcal{H}_2)(\log
R)^{d-1-j},
\label{eq:6.43}\end{split}\end{equation}
where
\begin{equation}
\label{eq:6.44}
\mathcal{B}_j(s_1,\mathcal{H}_1,\mathcal{H}_2)
= 
{d-1\choose j} \sum^j_{\nu = 0} {j\choose \nu}
\frac{\partial^{j-\nu}}{\partial s^{j-\nu}_2} D(s_1,
s_2)\Big|_{s_2 = -s_1}  \cdot
\frac{(-1)^{\nu} (v+1) \dots (v+\nu)}{(-1)^{\nu + v + 1} s^{u+v+\nu+2}_1}.
\end{equation}
We thus obtain
\begin{equation}
I_2 = \frac{1}{(d-1)!}
 \sum^{d-1}_{j=0}
\mathcal{C}_j(\mathcal{H}_1,\mathcal{H}_2)(\log R)^{d -1 -j} +
O(e^{-c\sqrt{\log N}}),
\label{eq:6.45}
\end{equation}
where
\begin{equation}
\mathcal{C}_j(\mathcal{H}_1,\mathcal{H}_2) =\frac {1}{ 2\pi
i}\mathop{\int}_{\mathcal{L}_3}
\mathcal{B}_j(s_1,\mathcal{H}_1,\mathcal{H}_2)\, ds_1\quad
(j = 0,1,2,\dots, d - 1) .
\label{eq:6.46}
\end{equation}
It remains to estimate these quantities, which are independent
of $R$.

We are allowed to transform the contour $\mathcal L_3$ in
\eqref{eq:6.46} to the contour $\mathcal L'$, defined in
\eqref{eq:4.4}.
Our task
is now the estimation of the integral $\mathcal C_j$ on the new contour
$\mathcal L' = \mathcal L'_0 \cup \mathcal L'_1$
since the integral on the horizontal segments $|t| = U$ is $O
\bigl(e^{-c \sqrt{\log N}}\bigr)$.

\section{Comparison of $D(s, -s)$ and $D(0,0)$}
\label{sec:D}

We have seen in Section~\ref{sec:A} that
by Lemma~\ref{l:A1} we can estimate
$\frac{\partial^i}{\partial s^i_1} \frac{\partial^j}{\partial
s^j_2} D(s_1, s_2)$ with the aid of $D(s_1, s_2)$.
We will show now how to estimate $|D(s, - s) | / D(0,0)$ from
above when $s$ is on the contour $\mathcal L'$.
This, together with Lemma~\ref{l:D1} will play a crucial role in the
estimation of $I_{2,1}$ which is the main part of~$I_2$.

First we note that if $s \in \mathcal L'$ is on the semicircle
$\mathcal L'_0$,
then by \eqref{eq:A12} we obtain
\begin{equation}
|D(s, -s) | \leq e^{C\sqrt K} D(0, 0), \qquad (s \in \mathcal L'_0).
\label{eq:D1}
\end{equation}

Thus, in the following we may suppose
\begin{equation}
s = it, \qquad t > 0,
\label{eq:D2}
\end{equation}
since $|D(-it, it)| = |D(it, - it)|$.

First we will examine  the behavior of the functions $D_0(s,
-s)$ and $G_1(s, -s)$ on the imaginary axis, which requires a lemma concerning $W(s)$ from \eqref{eq:74}. 
\begin{lemma}
\label{l:D1}
There exist positive absolute constants $t_0$ and $t_1 > 1$ such that
\begin{equation}
|W(it)| \geq e^{t^2/6} \geq
1 = W(0) \quad \text{for } |t| \leq t_0,
\label{eq:D3}
\end{equation}
\begin{equation}
|W(it)| \geq t^{2/3} \qquad\quad \text{for } |t| \geq t_1.
\label{eq:D4}
\end{equation}
\end{lemma}

\noindent
{\it Proof.}
We will use that in a neighborhood of $s = 0$ we have for the
entire function $W(s)$ the representation
\begin{equation}
W(s) = 1 + \gamma_0 s + \sum^\infty_{\nu = 1} \gamma_\nu s^{\nu
+ 1}
\label{eq:D5}
\end{equation}
where $\gamma_0 = \gamma$ is Euler's constant and (see \cite{Ivic},
Notes on p.~49)
\begin{equation}
\gamma_0 = \gamma = 0.5772157\dots, \quad \gamma_1 =
0.07281\dots\ .
\label{eq:D6}
\end{equation}
This implies
\begin{align}
|W(it)|^2 &= W(it)W(-it) = (1 + i\gamma t - \gamma_1 t^2 + O(t^3))
(1 - i\gamma t - \gamma_1 t^2 + O(t^3))
\label{eq:D7}\\
&= 1 + t^2(\gamma^2 - 2\gamma_1) + O(t^3) \quad \text{if } t \to 0.
\nonumber
\end{align}
Now \eqref{eq:D6}--\eqref{eq:D7} prove \eqref{eq:D3} for $|t| \leq t_0$, while
\eqref{eq:D4}
clearly holds by \eqref{eq:4.2}.

\begin{rema}
If \eqref{eq:D3} is true for any $t$ (which could be checked by
computers, since $t_0, t_1$ are explicitly calculable), then the
following simple lemma is not necessary.
\end{rema}

\begin{lemma}
\label{l:D2}
Given any positive constants $B_0, B_1, \varepsilon $ we have for any
$t \in [B_0, B_1]$ and any $X > C(B_0, B_1, \varepsilon )$
\begin{equation}
J(t,X) := \prod_{p \leq X} \frac{|1 - p^{-1 - it}|}{1 - p^{-1}}
\geq c(B_0, B_1) (\log X)^{1/2 - \varepsilon }.
\label{eq:D8}
\end{equation}
\end{lemma}

\noindent
{\it Proof.}
Let us fix $t$.
Since every factor is at least $1$, we can neglect those with
$\cos(t \log p) > 0$.
On the other hand, if $\cos(t \log p) \leq 0$, then we have
\begin{equation}
\log \left( \frac{|1 - p^{-1 - it}|}{1 - p^{-1}} \right) > \log
\frac1{1 - p^{-1}} > \frac1p.
\label{eq:D9}
\end{equation}
The primes satisfying $\cos (t\log p) \leq 0$ are in intervals
of type
\begin{equation}
I_j = \left[ \exp \left( \frac{\pi\left(2j +
\frac12\right)}{t}\right), \exp \left(\frac{\pi\left(2j +
\frac32\right)}{t}\right) \right] =: [e^{m_j}, e^{m_j + \pi / t}]
\label{eq:D10}
\end{equation}
and $I_j \subset [1, X]$ if $2\pi(j+ 3/4)/t \leq \log X$, that
is, if
\begin{equation}
j \leq \frac{t \log X}{2\pi} - \frac34 =: j^*.
\label{eq:D11}
\end{equation}

Using the prime number theorem we obtain by partial summation
\begin{equation}
\sum_{p \in I_j} \frac1p \sim \int\limits^{e^{m_j + \pi /
t}}_{e^{m_j}} \frac{dx}{x \log x} = \log \frac{2j + \frac32}{2j +
\frac12} = \frac1{2j} + O \left(\frac1{j^2}\right).
\label{eq:D12}
\end{equation}
Hence, by \eqref{eq:D9} we have
\begin{equation}
\log J(t,X) > \sum_{1 \leq j \leq j^*} \frac{1 - \varepsilon }{2j}
+ O(1) > \frac{1 - \varepsilon }{2} \log_2 X - c'(B_0, B_1).
\label{eq:D13}
\end{equation}

\begin{rema}
Working more carefully we could prove Lemma~\ref{l:D2} with $(\log
X)^{\frac12 - \varepsilon }$ replaced by $\log X$. But actually any
lower bound larger than
$C(t_0, t_1) = \max\limits_{t_0 \leq t \leq t_1} |W(it)|^{-1}$
would suffice for us.
\end{rema}

Taking into account the trivial relation
\begin{equation}
\big|1 - p^{-1 - it} \big|^{-1} \leq \big|1 - p^{-1}\big|^{-1},
\label{eq:D14}
\end{equation}
we obtain from Lemmas \ref{l:D1}, \ref{l:D2} the following

\begin{lemma}
\label{l:D3}
We have, with a sufficiently small constant $t_0 < 1$ and suitable positive
constants $c'$ and $c''$,
\begin{equation}
E_0(t) := \left| \frac{D_0(it, -it) G_1(it, -it)}{D_0(0,0)
G_1(0,0)} \right| \leq e^{-c'(a + b)t^2}, \quad \text{if } |t|
\leq t_0,
\label{eq:D15}
\end{equation}
and for any $t > t_0$ and $N > N_0$,
\begin{equation}
E_0(t) \leq e^{-c''(a + b)} \max(1, |t|)^{-(a + b)/2}.
\label{eq:D16}
\end{equation}
\end{lemma}

\noindent{\it Proof.}
By \eqref{eq:D14} and the definition of
$D_0$ in \eqref{eq:74} we clearly have
\begin{equation}
|G_1(it, -it)| \leq G_1(0,0),
\label{eq:10.16a}
\end{equation}
\begin{equation}
|D_0(it, -it)| = |W(it)|^{-(a + b)}, \quad D_0(0,0) = W(0) = 1,
\label{eq:10.16b}
\end{equation}
which immediately imply
\begin{equation}
E_0(t) \leq |W(it)|^{-(a + b)} .
\label{eq:10.16c}
\end{equation}
Hence, by Lemma~\ref{l:D1} we have \eqref{eq:D15} and
\eqref{eq:D16} for $|t| \geq t_1$.
Finally, for $t_0 \leq |t| \leq t_1$ we have by Lemma~\ref{l:D2}
\begin{equation}\begin{split}
E_0(t) &= \bigl(J(|t|, V) |W(it)| \bigr)^{-(a + b)}
\leq \bigl(c(t_0, t_1) ( \log V)^{1/3} \cdot C^{-1} (t_0, t_1)
\bigr)^{-(a + b)} \\&
\leq (c \log_2 N)^{-(a + b)/3} 
\leq (e^{-c''}t_1)^{-(a + b)},
\label{eq:10.16d} \end{split}\end{equation}
which proves \eqref{eq:D16}.

We will continue our study of $D(it, -it)$ with that of
\begin{equation}
L_4(t) := \log \frac{|G_4(it, -it)|}{G_4(0,0)} = \text{\rm
Re}\log \frac{G_4(it, -it)}{G_4(0,0)}.
\label{eq:D17}
\end{equation}

We first divide each term by $\left(1 + \frac{\nu_3 (p)}{p}\right)
\left(1 - \frac{1}{p}\right)^d$ in the product representation of
both $G_4(it, -it)$ and $G_4(0,0)$. After this we take the
logarithm of each term and use the formula 
\begin{equation}
\log(1 - z) = - \bigg( z + \sum^\infty_{m = 2} \frac{z^m}{m} \bigg),
\ \text{ if }\ |z| < 1.
\label{eq:D18}
\end{equation}
Now we separate the effect of the linear terms and those of
order $m \geq 2$ and write accordingly
\begin{equation}
L_4(t) = L_{4,1}(t) + L_{4,2}(t).
\label{eq:D19}
\end{equation}

We have by the trivial relations $\nu_1(p) \leq a$, $\nu_2(p)
\leq b$,
\begin{equation}
L_{4,1}(t) = \sum_{p > V} \left(\frac{a+b}{p} - \frac{\nu_1(p) +
\nu_2(p)}{p(1 + \nu_3(p)/p)}\right) (\cos(t\log p) - 1) \leq 0.
\label{eq:D20}
\end{equation}
(We remark that the sum is convergent, since $\nu_1(p) = a$,
$\nu_2(b) = p$ for $p\nmid \Delta$.)

The logarithms of the higher order terms of $G_4$ which do not involve
the functions $\nu_i(p)$
 can be estimated from above in modulus for any $t$
by
\begin{equation}
(a + b) \sum_{p > V} \sum^\infty_{m = 2} \frac2{mp^m} \leq C(a +
b) \frac1{V \log V} \leq \frac{C}{(\log_2 N)^3}.
\label{eq:D21}
\end{equation}

Similarly we have for the contribution of the numerator to
$L_{4,2}(t)$ the upper estimate (valid for any $t$)
\begin{equation}
\sum_{p > V} \sum^\infty_{m = 2} \frac{2(a + b)^m}{mp^m} \leq
\frac{C(a + b)^2}{V \log V} \leq \frac{C(a + b)}{(\log_2 N)^3} .
\label{eq:D22}
\end{equation}

We see from \eqref{eq:D16} that this estimation is
sufficient for $|t| \geq t_0$ but not for small values of $t$.
Then, working more carefully we have for the contribution of the
terms of $G_4$
involving the functions $\nu_i(p)$
to $L_{4,2}(t)$ the upper estimate
\begin{align}
&\sum_{p > V} \sum^\infty_{m = 2} \sum^m_{j = 0} {m \choose j}
\frac{(\nu_1(p))^j (\nu_2(p))^{m - j}}{m(p + \nu_3(p))^m}
\bigl(1 - \cos((m - 2j) t \log p)\bigr)
\label{eq:D23}\\
&\leq C \sum_{p > V} \sum^\infty_{m = 2} \frac{(a + b)^m}{m p^m}
t^2 m^2 \log^2 p \leq C(a + b)^2 t^2 \sum_{p > V} \frac{\log^2
p}{p^2} \nonumber\\
&\leq \frac{C(a + b)^2 t^2 \log V}{V} \leq \frac{C(a +
b)t^2}{\log_2 N}. \nonumber
\end{align}

It is easier to see that 
the contribution of the
terms of $G_4$ which do not involve the functions $\nu_i(p)$
to $L_{4,2}$ are majorized by 
\begin{equation}
C(a + b)t^2 \sum_{p > V} \frac{\log^2 p}{p^2} \leq \frac{Ct^2}{\log_2N}.
\label{eq:D24}
\end{equation}
Summarizing \eqref{eq:D17}--\eqref{eq:D24} we have proved

\begin{lemma}
\label{l:D4}
$\frac{|G_4(it, -it)|}{G_4(0,0)} \leq \exp \left(C\frac{(a +
b)}{\log_2 N} \min(1, t^2)\right)$.
\end{lemma}

Comparing the above with \eqref{eq:D15}--\eqref{eq:D16} we see
that \eqref{eq:D15}--\eqref{eq:D16} remain valid if we multiply
them by $G_4(it, -it)/G(0,0)$.
This proves the final result of this section:

\begin{lemma}
\label{l:D5}
$|D(it, -it)| \leq \max(1,|t|)^{-(a + b)/2} D(0,0)$ for any real
$t$.
\end{lemma}

Together with \eqref{eq:D1} this implies

\begin{lemma}
\label{l:D6}
$|D(s, - s)| \leq e^{C\sqrt K} \max(1, |t|)^{-(a + b)/2} D(0,0)$
for $s \in \mathcal L'$.
\end{lemma}

\section{Estimate of $I_2$. Evaluation of $I$}
\label{sec:E}

In this section we will estimate the integral $I_{2,1}$ based on
formulas \eqref{eq:6.44}--\eqref{eq:6.46},
using Lemmas~\ref{l:A1} and \ref{l:D6}.

First we obtain from the above lemmas for $s \in \mathcal L'$ by
$j \leq d \leq \min(a,b) \leq K$, $v \leq \sqrt K$
\begin{align}
\mathcal B_j(s, \mathcal H_1, \mathcal H_2)
&\ll \frac{d^j}{|s|^{u + v + 2}} \sum^j_{\nu = 0} (CK(\log_2 N +
\log T))^{j - \nu} \prod^\nu_{i = 1} \frac{\left(\frac{v}{i} +
1\right)}{|s|} |D(s, - s)|
\label{eq:E1}\\
&\ll \frac{e^{C\sqrt K} D(0,0) d^j(\log(|t| + 3))^j}{\max(1,
\sqrt{|t|})^{a + b}} \frac{\delta^{-(u + v)}_0}{|s|^2} \sum^j_{\nu
= 0} (CK \log_2 N)^{j - \nu} (K \log_2 N)^\nu \nonumber\\
&\ll e^{C\sqrt K} (CK^2 \log_2 N)^j \cdot
\frac{\delta^{-(u+v)}_0}{|s|^2} D(0,0).
\nonumber
\end{align}

Integrating the above upper bound along $\mathcal L'$ we obtain
\begin{equation}
\mathcal C_j(\mathcal H_1, \mathcal H_2) \ll e^{C\sqrt K}(CK^2
\log_2 N)^j \delta^{-(u+v+1)}_0 D(0,0).
\label{eq:E2}
\end{equation}

Finally, summation over $j \leq d - 1$ yields in \eqref{eq:6.45}
by $R \gg N^c$
\begin{align}
I_{2,1} &\ll \frac{e^{C\sqrt K} D(0,0) \delta^{-(u+v+1)}_0 (\log
R)^{d - 1}}{(d - 1)!} \cdot \sum^{d - 1}_{j = 0} \left(
\frac{CK^2 \log_2 N}{\log R}\right)^j
\label{eq:E3}\\
&\ll \frac{e^{C(u + v)} D(0,0) (\log R)^{d - 1} (\sqrt K\log_2
N)^{u + v + 1}}{(d - 1)!} \nonumber\\
&\ll \frac{D(0,0)(\log R)^{d + u + v}}{(d + u + v)!} \cdot
\left( \frac{CK^{3/2} \log_2 N}{\log R} \right)^{u + v + 1}\nonumber\\
&\ll \frac{D(0,0)(\log R)^{d + u + v}}{(d + u + v)!} (\log
N)^{-\sqrt K /50}. \nonumber
\end{align}

This implies by \eqref{eq:6.42} and \eqref{eq:6.45}
\begin{equation}
I_2 \ll \frac{D(0,0)(\log R)^{d + u + v}}{(d + u + v)!} (\log
N)^{-\sqrt K / 50} + e^{-c\sqrt{\log N}}.
\label{eq:E4}
\end{equation}

This yields by Lemma~\ref{l:B2} the final asymptotic evaluation of $I$
in \eqref{eq:73} by $D(0,0) = G(0,0)$
as
\begin{equation}
I = \frac{{v + u\choose u} (\log R)^{d + v + u}G(0,0)}{(d + v +
u)!} \left(1 + O\left( \frac{K\bar d^* \log_2 N}{\log
R}\right)\right) + O \big(e^{-c\sqrt{\log N}}\big),
\label{eq:E5}
\end{equation}
where
\begin{equation}
G(0,0) = \prod_{p \nmid P} \left(1 - \frac{\nu_p(\mathcal
H)}{p}\right) \prod_p \left(1 - \frac1p\right)^{-|\mathcal H|} =
\frac{{\mathfrak S}(\mathcal H) P}{|A(\mathcal H)|},
\label{eq:mar11.6}
\end{equation}
thereby proving Theorem~\ref{ujth:1}.

\section{A Bombieri--Vinogradov type theorem}
\label{sec:11}

In the present section we will prove a
modified Bombieri--Vinogradov theorem,
where the examined moduli are all multiples of a single modulus~$M$.
It would facilitate our task if we were entitled to use
the following hypothesis.

\begin{HypoSY}
If $L(1 - \delta, \chi) = 0$ for a $\delta > 0$ and a real
primitive character $\chi(\text{\rm mod}\, q)$, $q \leq Y$, then
\begin{equation}
\delta \geq \frac1{3\log Y}
\label{eq:11.1}
\end{equation}
for $Y > C_0$, an explicitly calculable absolute constant.
\end{HypoSY}

We note that we have the effective unconditional estimate (\cite{GS},
\cite{Pi1}), valid for $q > q_0$:
\begin{equation}
\delta \geq \frac{1}{\sqrt q}.
\label{uj8.2}
\end{equation}

A further observation (similar to that of Maier \cite{Ma}) is that
by the Landau--Page theorem (cf.\ Davenport \cite[\S 14]{Da}),  with some constant
$c$ in place of $1/3$, or Pintz \cite{Pi2} with $(1/2 + o(1))$)
for any given $Y$ there is at most one real primitive character
$\chi_1$ which does not fulfill \eqref{eq:11.1}.
This makes it possible to turn Hypothesis~$S(Y)$ into a theorem,
valid for a sequence $Y = Y_n \to \infty$ (for $n > n_0$, an
explicitly calculable absolute constant) with
\begin{equation}
Y_n \leq \exp \left(\sqrt{Y_{n-1}}\right).
\label{uj8.3}
\end{equation}

In order to show this, suppose that \eqref{eq:11.1} is false for a
sufficiently large $Y'$, i.e.\ by \eqref{uj8.2} there exists a
$\chi_1$ $\text{\rm mod}\, q_1 \leq Y'$ such that $ L(1 -
\delta_1, \chi_1) = 0$ with
\begin{equation}
\frac1{\sqrt{Y'}}  \leq \min \left(
\frac1{\sqrt{q_1}}, c_0 \right) \leq \delta_1 < \frac1{3\log
Y'} .
\label{uj8.4}
\end{equation}
Let us choose $\widetilde Y > Y'$ in such a way, that
\begin{equation}
\widetilde Y = \exp \left(\frac1{3\delta_1}\right)
\Leftrightarrow
\delta_1 = \frac1{3\log\widetilde Y}.
\label{uj8.5}
\end{equation}
Then for any other zero $1 - \delta_2$ belonging to a real
primitive $\chi_2$ $\text{\rm mod}\, q_2$, $q_2 \leq \widetilde Y$,
we have by the Landau--Page theorem in the version of Pintz \cite{Pi2}
\begin{equation}
\max(\delta_1, \delta_2) > \frac1{3\log \widetilde Y}
\Leftrightarrow \delta_2 > \frac{1}{3\log \widetilde Y}.
\label{uj8.6}
\end{equation}

Now, \eqref{uj8.4}--\eqref{uj8.6} show that \eqref{eq:11.1} is true for
a value $Y = \widetilde Y$ satisfying
\begin{equation}
Y'< \widetilde Y < \exp \bigl(\sqrt{Y'}/3\bigr).
\label{uj8.7}
\end{equation}
We can formulate this as

\begin{lemma}
\label{l:11.1}
Hypothesis $S(Y)$ holds for a sequence $Y_n \to \infty$ with
\begin{equation}
Y_n \leq \exp \big(\sqrt{Y_{n - 1}}\big)
\label{eq:11.2}
\end{equation}
where $Y_0$ can be chosen with $Y_0 < C_0$, an explicitly
calculable absolute constant.
\end{lemma}

An alternative to this Lemma and this approach would be to use Heath-Brown's theorem \cite{HB-Siegel}
(but only in case of Theorem~1)
according to which either

(i)~$S(Y)$ holds for every $Y > C$, with some absolute
constant $C$,

\noindent
or

(ii)~there are infinitely many twin primes.

\noindent The significance of the real zeros in Hypothesis $S(Y)$ is that a
similar inequality holds with $\text{\rm Re}\, \varrho$ in place
of $1 - \delta$ and with a constant $c_0$ in place of $1/3$ if
$\text{\rm Im}\,\varrho$ is not too large; 
this is the standard zero-free region of $L$-functions
(cf.\ Davenport \cite[\S 14]{Da}).

\begin{lemma}
\label{l:11.2}
There exists an explicitly calculable absolute constant
$c_0 < 1/3$
such that $L(s, \chi ) \neq 0$ in the region
\begin{equation}
\sigma > 1 - \frac{c_0}{\log(q(|t| + 3))},
\label{eq:11.3}
\end{equation}
apart from possible real exceptional zeros of real $L$-functions.
\end{lemma}

Now we are in a position to formulate and prove the following
theorem (which is similar to but stronger than Lemma~6 of Maier
\cite{Ma}).

\begin{theorem}
\label{th:3}
Let $c^*$ be an arbitrary, fixed constant.
Let $Y = Y(X)$ be a strictly monotonically increasing function
of $X$ with
\begin{equation}
\exp\bigl(2 \sqrt{\log X}\bigr) \leq Y(X) \leq X.
\label{eq:mar12.10}
\end{equation}
Then there exists a sequence $X_n \to \infty$ satisfying
$X_1 < C'_0$, an explicitly calculable constant,
with the following property.
Let $X = X_n$, $\mathscr L = \log X$, $M$ be a natural number
$\leq \min\bigl(\sqrt{Y(X)}/4, X^{1/8}\bigr)$,
\begin{equation}
Q^* = X^{1/2} M^{-3} \exp \big( - c^* \sqrt{\log X}\big),
\label{eq:11.5}
\end{equation}
\begin{equation}
E^*(X, q) := \max_{x \leq X} \max_{(a,q) = 1}
|E(X, q, a)| :=  \max_{x \leq X} \max_{(a,q) = 1} \Big|
\sum_{\substack{p \leq x\\ p \equiv a(\text{\rm
mod}\,q)}} \log p - \frac{x}{\varphi(q)}\Big| ,
\label{eq:11.6}
\end{equation}

Then
\begin{equation}
\sum_{\substack{q \leq Q^*\\ (q,M) = 1}} E^*(X, Mq) \ll 
\frac{X}{M} \mathscr L^{15} \exp \biggl(-c_2 \frac{\log X}{\log Y(X)} \biggr),
\label{eq:11.8}
\end{equation}
where $c_2 = \min(c^*/6, c_0/4)$.
\end{theorem}

\begin{rema}
The above theorem holds with any $X$, for which
$S(Y(X)) = S(Y)$ is true, i.e.
\begin{equation}
L(s, \chi ) \neq 0 \quad\text{ for } \
s \in \left( 1 - \frac1{3\log Y}, 1 \right]
\label{eq:11.8a}
\end{equation}
holds without exception for all real primitive characters $\chi$
$\text{\rm mod}\, q$, where
\begin{equation}
q \leq Y = Y(X).
\label{eq:11.8b}
\end{equation}
\end{rema}

\noindent
{\it Proof.}
We will choose our sequence $X_n = Y^{-1}(Y_n)$, where $Y_n$ is
the sequence supplied by Lemma~14 (for which $S(Y)$ is true)
and $Y^{-1}$ is the inverse function of $Y(X)$.
Alternatively, if \eqref{eq:11.8a}--\eqref{eq:11.8b} hold, then
we can choose $X$ as an arbitrary sufficiently large number.
In both cases $S(Y)$, i.e.\ \eqref{eq:11.8a}--\eqref{eq:11.8b} hold.
Using the explicit formula for primes in arithmetic progressions
with $T^* = \sqrt X \log^2 X$ $(\varrho = \beta + i \gamma = 1
- \delta + i \gamma$ denotes a generic zero of an $
L$-function) we obtain (cf.\ Davenport \cite[\S 19]{Da} for any
$a$ with $(a,q) = 1$, $q \leq Q^*$, $y \leq X$ the relation
\begin{equation}
E(y, q, a) = -\frac1{\varphi(q)}
\sum_{\chi (q)} \overline{\chi} (a)
\sum_{\substack{\varrho = \varrho_\chi\\
\beta \geq 1/2, |\gamma|
\leq T^*}} \frac{y^\varrho}{\varrho} + O\big(\mathcal L^2\sqrt{y}\big) .
\label{7.2}
\end{equation}
The effect of the last error term is clearly suitable,
$O(Q^* \mathcal L^2 \sqrt X)$ in total.
We can classify zeros of all primitive $  L$-functions
$\text{\rm mod}\, q\widetilde M \leq Q^*M, (\widetilde M \mid M)$, 
up to height $T^*$ into
$O(\mathscr L^4)$ classes
$B(\kappa, \lambda, \mu, \nu)$ by Lemma~15, as
\begin{equation}
\widetilde M \in [M_\lambda/2, M_\lambda),\quad
q \in [Q_\nu / 2, Q_\nu),\quad
\gamma \in [T_\mu / 2, T_\mu), \quad
\delta \in \left[ \frac{\kappa c_0}{\mathscr L},
\frac{(\kappa + 1)c_0}{\mathscr L}\right),
\label{7.3}
\end{equation}
where
\begin{equation}
M_\lambda = 2^\lambda \leq 2M, \quad
Q_\nu = 2^\nu \leq 2 Q^*, \quad
T_\mu = 2^\mu \leq 2T^*, \quad
\frac{\kappa c_0}{\mathscr L} \leq \frac12,
\label{7.4}
\end{equation}
with the additional class of index $0$: $\gamma \in [0,1) = [0, T_0)$.
The set of quadruples $\kappa, \lambda, \mu, \nu$ satisfying
\eqref{7.4} with $\nu \geq
1$, $\mu \geq 0$, $\lambda \geq 1$, $\kappa \geq 0$ will be
 denoted by~$\mathcal B$.

In this case we have clearly by \eqref{7.2}, similarly to
Davenport~\cite[\S 28]{Da},
\begin{equation}
\sum_{\substack{q \leq Q^*\\
(q, M) = 1}} E^*(X, qM)
\ll \frac{X}{M} \mathscr L^6 \max_{\kappa, \lambda,
\mu, \nu \in \mathcal B}
\frac{N^*(1 - \frac{(\kappa + 1) c_0}{\mathscr L}, M_\lambda Q_\nu,
T_\mu)}{Q_\nu T_\mu} X^{- c_0\kappa/\mathscr L},
\label{7.5}
\end{equation}
where
\begin{equation}
N^* (\sigma, Q, T) = \sum_{Q/2 < q \leq Q}
\sum_{\substack{\chi (q) \\ \chi \text{ primitive}}}
\sum_{\substack{\varrho = \varrho_{\chi}\\ \beta \geq \sigma, |\gamma |
\leq T}} 1.
\label{uj7.6}
\end{equation}
We will see that,
in order to prove our theorem, it will be enough to prove for any
quadruple
$\delta, M'Q, T$ with the property
(cf.\ \eqref{eq:11.3}, \eqref{eq:11.8a}--\eqref{eq:11.8b})

\begin{equation}
\gathered
\frac{c_0 }{\log (Q M' T)} \leq \delta \leq \frac12, \quad
2 \leq M' \leq M, \quad
2 \leq Q \leq Q^*, \quad
1 \leq T \leq T^* \ \text{ or}\\
\frac{c_0}{\log Y} \leq \delta \leq \frac12, \quad
2 \leq M' \leq M, \quad
Q \leq \sqrt{Y}, \quad
T = T_0 = 1,\\
0 \leq \delta \leq \frac12, \quad
2 \leq M' \leq M, \quad
Q > \sqrt{Y}, \quad T = T_0 = 1
\endgathered
\label{7.6}
\end{equation}
the crucial inequality
\begin{equation}
N^*(1 - \delta, M' Q, T) \ll
\mathcal L^9 Q T X^\delta \exp(-c \log X/ \log Y)
\label{7.7}
\end{equation}
with some positive absolute constant $c$.
The first line in \eqref{7.6} is meant to cover all non-real
zeros, the second
and third lines are meant to cover the real zeros.

We will use Theorem~12.2 of
Montgomery \cite{Mo}
\begin{equation}
N^*(1 - \delta, Q, T) \ll (Q^2 T)^{\frac{3\delta}{1 + \delta}}
(\log QT)^9.
\label{7.8}
\end{equation}
(We do not need
for the range $\delta \leq 1/5$
the stronger inequality of Theorem 12.2 of \cite{Mo}
with the exponent $3\delta/(1+\delta)$ replaced by
the smaller
 $2\delta/(1 - \delta)$.)
Since $3\delta/(1 + \delta) \leq 1$, \eqref{7.7} will follow if
we can show
\begin{equation}
M^{6\delta}(Q)^{\frac{6\delta}{1+\delta} - 1} \ll X^\delta e^{-c_2
\sqrt{\log X}} \text{ with } c_2 =  c^*/6.
\label{7.9}
\end{equation}
Since in the range $0 \leq \delta \leq 1/2$ we have
$\frac{6\delta}{1 + \delta} - 1 \leq 2\delta$, this is
true by the definition $Q^* = X^{1/2} M^{-3} \exp ({-c^* \sqrt{\log
X}})$, if $\delta \geq 1/12$.

In case of $\delta \leq 1/12$ we have by \eqref{7.8}
\begin{equation}
N^*(1 - \delta, M' Q, T) \ll (Q T)^{1/2} M^{6\delta}.
\label{7.10}
\end{equation}

If we have here $QT \geq \exp(\sqrt{\log X})$, then \eqref{7.10}
directly implies \eqref{7.7}, since
\begin{equation}
\frac{N^*(1 - \delta, M' Q, T)}{QT} \ll (QT)^{-1/2}
M^{6\delta}
\ll  X^\delta \exp
\big( - \sqrt{\log X}/2\big).
\label{7.11}
\end{equation}

If $\delta \leq 1/12$ and
$QT \leq e^{\sqrt{\log X}} \leq \sqrt{Y}$, then
$\delta \geq c_0 / \log (MQT)$ or $\delta \geq \frac1{3\log Y} >
\frac{c_0}{\log Y}$ by \eqref{eq:11.8a}--\eqref{eq:11.8b},
since the modulus of the corresponding primitive character is $q
\widetilde M \leq 4QM \leq Y$.
Hence,
\begin{align}
\left(\frac{M^6}{X}\right)^\delta
&\leq X^{-\delta / 4} \leq \exp
\left( - \frac{c_0}{8} \min \left( \frac{\log X}{\log QT} ,
\frac{\log X}{\log M} \,, \frac{\log X}{\log \sqrt Y}\right)\right)
\label{eq:11.13}\\
& =\exp \left( -\frac{c_0}{8} \frac{\log X}{\log \sqrt{Y}}\right). 
\nonumber
\end{align}

\begin{rema}
The condition for $Q^*$ could be weakened to $Q^* < X^{1/2}
M^{-1} \exp(-cf(X,Y))$ but this has no significance in our application.
\end{rema}

\section{Proof of Theorem \ref{ujth:2}}
\label{sec:12}

The method of proof of Theorem~\ref{ujth:2} is quite similar to that of
Theorem~\ref{ujth:1}. The basic difference is that instead of the trivial
problem of the distribution of integers in arithmetic progressions
we have to use properties of the distribution of primes in
arithmetic progressions. Since we have to consider the (weighted)
sum of the error terms in the formula for the number of primes in
arithmetic progressions, the Bombieri--Vinogradov theorem can help
us. However,
due to the relatively weak estimate of
the original Bombieri--Vinogradov theorem, it does not lead to better
results than $\liminf\limits_{n \to\infty} (p_{n+1} - p_n) / \log p_n = 0$.
That is partly
why we need  to use  Theorem~6 instead.
Our situation is even more complicated here, since we need
 the moduli of the progressions to be
the multiples of a number~$V$. Fortunately our present Theorem~\ref{th:3}
solves this problem in a completely satisfactory way, even without
loss if $P = M \leq \exp((1 + o(1)) \sqrt{\log N})$
which is now the case by $V = \sqrt{\log N}$.

We will suppose that $N = X_n / 3$, $n$ is sufficiently large,
and $M = P$ in Theorem~\ref{th:3}.
(If we use Heath-Brown's theorem \cite{HB-Siegel} we may assume
Hypothesis~$S(Y)$ for any $N$
and then $N$ can be an arbitrary, sufficiently large integer.)

In the course of proof we will follow closely the analogous
proofs of Propositions~4 and~5 in Sections 7--9 of
\cite{GPY}, so we will sometimes omit details.
Let
\begin{align}
\boldsymbol\Theta (x;q,a) :&= \sum_{\substack{p\le x\\ p\equiv
a (\textrm{mod}\ q) }}\log p  = [(a,q)=1]\frac{x}{\phi(q)} +
E(x;q,a),
\label{eq:12.1}
\end{align}
where $[S]$ is $1$ if the statement $S$ is true and $0$ if $S$
is false. We have
for a regular residue class $\widetilde a$ with respect to
$\mathcal H$ and $P$
\begin{align}
&\widetilde S_R (N; \mathcal H_1, \mathcal H_2, \ell_1, \ell_2,
P, \widetilde a,
h_0) := \sum^{2N}_{n = N + 1} \Lambda_R(n; \mathcal H_1, \ell_1)
\Lambda_R(n; \mathcal H_2, \ell_2) \theta (n + h_0)
\label{eq:12.2}\\
&=\! \frac{1}{(K\! +\! \ell_1)!(K\! +\! \ell_2)!}
\sum_{d,e \leq R} \!\! \mu(d) \mu(e)\!
\left(\!\log \frac{R}{d}\right)^{\!\! K + \ell_1}
\! \left(\!\log
\frac{R}{e}\right)^{\!\! K + \ell_2}
\!\!\!\!\!\!\!\! \sum_{\substack{1 \leq n \leq N, n
\equiv \widetilde a (\text{\rm mod}\, P)\\ d \mid P_{\mathcal H_1}(n),
e \mid P_{\mathcal H_2}(n)}}\!\!\!\!\!\!\! \theta(n + h_0). \nonumber
\end{align}

For the inner sum, we let $d=a_1a_{12}$, $e=a_2a_{12}$ where
$(d,e)=a_{12}$, and thus $a_1$, $a_2$, and $a_{12}$ are pairwise
relatively prime.
We may suppose in the following $(d, P) = (e, P) = 1$, otherwise
the last sum would be zero, since by the regularity of $\widetilde a$ we
have $(P_{\mathcal H_1}(n), P) = (P_{\mathcal H_2}(n), P) = 1$.
The $n$ for which
$d|P_{\mathcal{H}_1}(n)$ and  $e|P_{\mathcal{H}_2}(n)$ cover\marginpar{}
certain residue classes modulo $[d,e]$.
If  $n\equiv b'
(\text{mod} \, a_1a_2a_{12})$ is such a residue class, then
letting
$m=n+h_0 \equiv b'+h_0(\text{mod} \, a_1a_2a_{12})$,
$b \equiv b'(\text{\rm mod}\, a_1 a_2 a_{12})$, $b\equiv
\widetilde a(\text{\rm mod}\, P)$
we see this residue class contributes to the inner sum
\begin{equation}
\begin{split}
&\sum_{\substack{N + 1+h_0 \leq  m\le 2N+h_0\\ m\equiv b
+h_0\, (\text{mod} \, a_1a_2a_{12} P) }}
\!\!\!\!\!\!\!\!\! \theta(m)
= \theta(2N+h_0;a_1a_2a_{12} P, b+h_0) -
\theta(N + h_0;a_1a_2a_{12} P, b+h_0)\\
& =
[(b+h_0,a_1a_2a_{12}P) = 1]
\frac{N}{\phi(a_1a_2a_{12} P)} +
O\bigl(E^*(3N, a_1 a_2 a_{12} P)\bigr).\end{split}
\label{eq:7.7}
\end{equation}
We need to determine the number of these residue classes where
$(b+h_0,a_1a_2a_{12}P) = 1$ so that the main term is non-zero.
The condition $(\widetilde a + h_0, P) = (b + h_0, P) = 1$ is equivalent
to $\widetilde a$ being regular with respect to $\mathcal H^0$, since
$\widetilde a$ is regular with respect to $\mathcal H$.
Thus we will assume from now on that $\widetilde a$ is regular with respect
to~$\mathcal H^0$.
If $p|a_1$ then $b\equiv - h_j\ (\text{\rm mod}\, p)$ for some $h_j \in
\mathcal{H}_1$, and therefore
$b+h_0\equiv h_0 - h_j\ (\text{\rm mod}\,
p)$.  Thus, if $h_0$ is distinct modulo $p$ from all the $h_j\in
\mathcal{H}_1$ then all $\nu_p(\mathcal{H}_1)$ residue classes
satisfy the relatively prime condition, while otherwise
$h_0\equiv h_j(\text{\rm mod}\, p)$ for some $h_j\in \mathcal{H}_1$
leaving $\nu_p(\mathcal{H}_1)-1$ residue classes with a non-zero
main term. We introduce the notation
${\nu_p}^*({\mathcal{H}_1})$ for this number in either case,
where we define for a set $\mathcal{G}$
\begin{equation}
{\nu_p}^*(\mathcal{G}) =\nu_p(\mathcal{G}^0)-1.
\label{eq:7.8}
\end{equation}
and
\begin{equation}
\mathcal{G}^0 = \mathcal{G} \cup  \{h_0\}.
\label{eq:7.9}
\end{equation}
We extend this definition to
${\nu_d}^*({\mathcal{H}_1})$ for squarefree numbers $d$ by
multiplicativity. (The function ${\nu_d}^*$ is familiar in sieve
theory, see \cite{HR}.) The same applies for $\nu^*_d(\mathcal
H_2)$ and $\overline \nu^*_d \big((\mathcal H_1 \overline \cap
\mathcal H_2)\big)$, as in~\eqref{eq:5.2}.

Since
$E(n;q,a)\ll (\log N)$ if $(a,q)>1$ and $q \leq N$ we conclude
\begin{equation}
\begin{split}
&\sum_{\substack{N + 1\le n\le 2N, n \equiv
\widetilde a(P)\\
d|P_{\mathcal{H}_1}(n),  e|P_{\mathcal{H}_2}(n)}}
\theta(n+h_0) = {\nu_{a_1}}^*({\mathcal{H}_1})
{\nu_{a_2}}^*({\mathcal{H}_2}){\overline\nu_{a_{12}}}^*
\left((\mathcal{H}_1 \overline \cap
\mathcal{H}_2)\right)\frac{N}{\phi(a_1a_2a_{12}P)} \\
& \quad + O\left( d_K(a_1a_2a_{12})\left(
\big|E^*(3N ;a_1a_2a_{12} P)\big|
\right) \right).\end{split}
\label{eq:7.10}
\end{equation}

Let $\sum^{(P)}$ denote that the summation variables are
relatively prime to $P$ and to each other.
Substituting this into \eqref{eq:12.2} we conclude by $\ell_i \leq
K$
\begin{align}
&\tilde{\mathcal{S}}_R(N;\mathcal{H}_1,\mathcal{H}_2,
\ell_1, \ell_2, P, \widetilde a,
h_0)
\label{eq:7.11}\\
&
=\! \frac{N}{\varphi(P)(K\!\! +\! \ell_1)!(K\!\! +\! \ell_2)!} \!\!
\sum_{\substack{ a_1a_{12}\le
R\\ a_2a_{12}\le
R}}\!\!\!\!\!{}^{^{\scriptstyle (P)}}
\frac{\mu(a_1)\mu(a_2)\mu(a_{12})^2{\nu_{a_1}}^{\!\!\!*}({\mathcal{H}_1})
{\nu_{a_2}}^{\!\!\!*}({\mathcal{H}_2}){\overline
\nu_{a_{12}}}^{\!\!\!*}\left((\mathcal{H}_1 \overline \cap
\mathcal{H}_2)\right)}{\phi(a_1a_2a_{12})} \nonumber \\
&\hskip5cm
\times \left(\log \frac{R}{a_1a_{12}}\right)^{K + \ell_1}
\left(\log
\frac{R}{a_2a_{12}}\right)^{K + \ell_2} \nonumber \\
& \quad + O\left((\log R)^{4K}
\sum_{\substack{ a_1a_{12}\le R\\ a_2a_{12}\le R}}\!\!\!\!\!
{}^{^{\scriptstyle(P)}}\,
d_K(a_1a_2a_{12})
E^*(3N; a_1 a_2 a_{12} P)
\right) \nonumber \\
&
= \frac{N}{\varphi(P)} \tilde{ \mathcal{T}}_R(\mathcal{H}_1,\mathcal{H}_2,
\ell_1, \ell_2, h_0)
+O \left((\log R)^{4K}\mathcal{E}_K(N)\right).
 \nonumber
\end{align}

Using the notation $R^2 = Q^*$, we obtain from Theorem~\ref{th:3}
 by the trivial estimate
$|E(X, Pq, a)| \leq 2 q^{-1} P^{-1} X \log X$ (for $Pq \leq X$),
Lemma~\ref{lem:uj2a} and by H\"older's inequality with parameters
$\alpha = \nu + 1$, $\beta = (\nu + 1) /\nu$
where $\nu \in \mathbb Z^+$, $c' \log (K
+ 1) \leq \nu \leq c'' \log (K + 1)$,
 uniformly for $K \leq (\log N)/(2C)$,
\big($\sum^{\flat*}$ means summation
over squarefree integers which are relatively prime to~$P$\big)
\begin{align}
\allowdisplaybreaks
|\mathcal E_K(N)| &\leq \sum_{q \leq Q^*}{}^{^{\scriptstyle\flat*}}
d_K(q) E^* (3N
, Pq) \sum_{q = a_1 a_2 a_{12}} 1
\label{eq:8.17}\\
&\leq \sum_{q \leq Q^*}{}^{^{\scriptstyle\flat*}} d_K(q) d_3(q)E^*(3N, Pq)
\ = \ \sum_{q \leq Q^*}{}^{^{\scriptstyle\flat*}} \frac{d_{3K}(q)}{q^{1/\beta}}
\cdot q^{1/\beta } E^*(3N,Pq)\nonumber\\
&\leq \bigg(\sum_{q \leq Q^*}{}^{^{\scriptstyle\flat}}
\frac{(d_{3K}(q))^\beta}{q}\bigg)^{1/\beta}
\bigg(\sum_{q\leq Q^*}{}^{^{\scriptstyle\flat*}}
q^{\alpha/\beta} (E^*(3N, Pq))^\alpha\bigg)^{1/\alpha}\nonumber\\
&\leq \left(1 + \frac1 2 \log N\right)^{CK} (6N P^{-1} \log 3N)^{\nu /
(\nu + 1)} \bigg(\sum_{q\leq Q^*}{}^{^{\scriptstyle\flat*}} E^*(3N,
Pq)\bigg)^{\frac1{\nu +1}} \nonumber\\
& \ll (\log N)^{CK+1}
N P^{-1}\exp \left( - \frac{c_2 \sqrt{\log
N}}{\nu + 1} \right)\nonumber\\
&\leq N P^{-1} \exp\left( (CK + 1) \log_2 N - c_2(\nu + 1)^{-1} \sqrt{\log
N}\right)\nonumber \\&  \leq  N P^{-1}  \exp \left(-c\frac{ \sqrt{\log N}}{\log (K+1)}\right).
\nonumber
\end{align}
Since, by \eqref{eq:2.12}, $K$ satisfies the inequality
\begin{equation}
K \log_2 N < c \sqrt{\log N}/\log K .
\label{eq:8.18}
\end{equation}
From \eqref{eq:8.18} we have, finally
\begin{equation}
(\log R)^{4K} |\mathcal E_K(N)| \leq P^{-1} N \exp
\left(-c \frac{\sqrt{\log N} }{ \log (K+1)}\right).
\label{eq:8.19}
\end{equation}

So, our task is reduced to the evaluation of $\widetilde T_R$
which is very similar to $\mathcal T^*_R$ in \eqref{eq:73}.
Due to the more general treatment of $\mathcal T^*_R$ in
Section~5 than needed, the crucial part, the error analysis will
remain the same. The difference will be only the fact that we
have now $\varphi(a_1 a_2 a_{12})$ in the denominator in \eqref{eq:7.11}
in place of $a_1 a_2 a_{12}$.
Therefore $\frac{\nu_i(p)}{p^{1 + s_i}}$ has to be replaced by
$\frac{\nu_i(p)}{(p - 1)p^{s_i}}$ in the definition of $F(s_1,
s_2)$ and $G(s_1, s_2)$ in \eqref{eq:5.5} and
\eqref{eq:61} (where $s_i = s_1$,
$s_2$ or $s_3 = s_1 + s_2$).
However, factors of type $(1 - p^{-(1 + s_i)})$ remain
unchanged, since they arise from the zeta-factors.

Summarizing our results above we have
\begin{equation}
\widetilde S_R(N; \mathcal H_1, \mathcal H_2, \ell_1, \ell_2, P,
\widetilde a, h_0)
= \frac{N}{\varphi(P)} \mathcal T_R (\ell_1, \ell_2; \mathcal H_1,
\mathcal H_2) + O \left( \frac{N}{P} \exp \left( -
\frac{c\sqrt{\log N}}{\log_2 N}\right)\right)
\label{eq:12.3}
\end{equation}
where
\begin{equation}
\mathcal T_R(\ell_1, \ell_2, \mathcal H_1, \mathcal H_2) := \frac
1{(2\pi i)^2} \int\limits_{(1)}\!\int\limits_{(1)} F(s_1, s_2)
\frac{R^{s_1}}{s^{K + \ell_1 + 1}_1} \frac{R^{s_2}}{s^{K + \ell_2 +
1}_2}\, ds_1 ds_2,
\label{eq:12.4}
\end{equation}
\begin{equation}
F(s_1, s_2) := \prod_{p > V} \left(1 - \frac{\nu_1(p)}{(p -
1)p^{s_1}} - \frac{\nu_2(p)}{(p - 1)p^{s_2}} +
\frac{\nu_3(p)}{(p - 1)p^{s_3}} \right)
\label{eq:12.5}
\end{equation}
and for this paragraph we have with notation \eqref{eq:5.2} $(i =
1,2)$
\begin{equation}
\nu_i(p) = \nu^*_p(\mathcal H_i) = \nu_p(\mathcal H^0_i) - 1,
\quad \nu_3(p) = \bar \nu_p((\mathcal H_1 \bar\cap  \mathcal
H_2)^0) - 1 .
\label{eq:12.6}
\end{equation}

To factor out the dominant zeta-factors we write now, in place of
\eqref{eq:5.8}
\begin{equation}
F(s_1, s_2) = G_{\mathcal H_1, \mathcal H_2}(s_1, s_2)
\frac{\zeta(1 + s_1 + s_2)^{|(\mathcal H_1 \cap \mathcal H_2)^0|
- 1}}{\zeta(1 + s_1)^{|\mathcal H^0_1| - 1} \zeta(1 +
s_2)^{|\mathcal H^0_2| - 1}}
\label{eq:12.7}
\end{equation}
and define accordingly $G_i$ $(1 \leq i \leq 4)$ as in
\eqref{eq:5.9}--\eqref{eq:61} with
\begin{equation}
a = |\mathcal H^0_1| - 1, \quad b = |\mathcal H^0_2| - 1, \quad
d = |(\mathcal H_1 \cap \mathcal H_2)^0| - 1,
\label{eq:12.8}
\end{equation}
and with $\nu_i(p) p^{-s_i} / (p - 1)$ in place of
$\nu_i(p)p^{-1 - s_i}$.

Similarly to Section~9 of \cite{GPY} by symmetry we have to consider
three cases:

Case 1. $h_0 \notin \mathcal H \Longleftrightarrow a = K$, $b =
K$, $d = r$.

Case 2. $h_0 \in \mathcal H_1 \setminus \mathcal H_2
\Longleftrightarrow a = K - 1$, $b = K$, $d = r$.

Case 3. $h_0 \in \mathcal H_1 \cap \mathcal H_2
\Longleftrightarrow a = K - 1$, $b = K - 1$, $d = r - 1$.

\noindent (Cases 1 and 3 are basically the same.)

Since the results of the previous section are more general,
they apply to the error analysis here and we only have to evaluate
$G(0,0)$ in Cases 1--3.
Similarly to Section~9 of \cite{GPY} we have by \eqref{eq:12.6}
\begin{equation}
\nu_1(p) + \nu_2(p) - \nu_3(p) = \nu_p(\mathcal H^0_1) +
\nu_p(\mathcal H^0_2) - \bar \nu_p(\mathcal H^0_1 \bar \cap
\mathcal H^0_2) - 1 =
\label{eq:12.9}\\
 \nu_p(\mathcal H^0) - 1,
\end{equation}
\begin{equation}
a + b - d = |\mathcal H^0| - 1.
\label{eq:12.10}
\end{equation}
Hence, from the analogies of \eqref{eq:5.9}--\eqref{eq:61} we have now
\begin{equation}
G_1(0,0) = \prod_{p \leq V} \left(1 -
\frac{1}{p}\right)^{-(|\mathcal H^0| - 1)} =
\left(\frac{P}{\varphi(P)}\right)^{|\mathcal H^0| - 1},
\label{eq:12.11}
\end{equation}
\begin{align}
G_4(0,0) &= \prod_{p > V} \left(1 - \frac{\nu_p(\mathcal H^0) -
1}{p - 1} \right) \left( \frac{p}{p - 1} \right)^{|\mathcal H^0|
- 1} =
\label{eq:12.12}\\
&= \prod_{p > V} \left( \frac{p - \nu_p(\mathcal H^0)}{p}
\right) \cdot \left( 1 - \frac1{p} \right)^{-|\mathcal H^0|} : =
\bar{\mathfrak S}_V(\mathcal H^0).
\nonumber
\end{align}

Taking into account the term $\varphi(P)$ in the denominator in
\eqref{eq:12.3} we obtain
\begin{equation}
\frac{G(0,0)}{\varphi(P)} = \frac1{P} \prod_{p \leq V} \left(1 -
\frac1{p}\right)^{-|\mathcal H^0|} \bar {\mathfrak S}_V(\mathcal H^0).
\label{eq:12.14}
\end{equation}

Further we have by the comparison of
\eqref{eq:12.4},  \eqref{eq:12.7} and \eqref{eq:73}
\begin{equation}
u = K + \ell_1 - a = K + 1 - |\mathcal H^0_1| + \ell_1,
\quad v = K + 1 - |\mathcal H^0_2| + \ell_2, \quad
d = |\mathcal H^0_1 \cap \mathcal H^0_2| - 1.
\label{eq:12.15}
\end{equation}

The evaluation \eqref{eq:E5} of the crucial integral $I$ defined
in \eqref{eq:73} yields in our case \eqref{eq:12.3}--\eqref{eq:12.7}
the relation
\begin{align}
&\widetilde S_R(N; \mathcal H_1, \mathcal H_2, \ell_1, \ell_2, P,
\widetilde a, h_0)
=
\label{eq:12.16}\\
&= N  \frac{G(0,0)}{\varphi(P)} \frac{{v + u\choose u}
(\log R)^{d + v + u}}{(d + v + u)!} \left(1 + O \left(
\frac{K\bar d^* \log_2 N}{\log R} \right)\right)\nonumber\\
&\quad  + O
\left(\frac{N}{\varphi(P)} e^{-c \sqrt{\log N}}\right).
\nonumber
\end{align}

Let us observe that on the right-hand side the residue class
$\widetilde a$ does not appear at all.
Therefore we can add this together for all $|A(\mathcal H^0)|$
regular residue classes $\widetilde a(\text{\rm mod}\, P)$ with respect
to $\mathcal H^0$ and $P$, since the contribution of those with
$(\widetilde a +
h_0, P) > 1$ is zero, as mentioned after \eqref{eq:7.7}.
Taking into account the trivial relations \eqref{eq:2.17}--\eqref{eq:2.18} for
$\mathcal H^0$ in place of $\mathcal H$ we obtain from
\eqref{eq:12.14}
\begin{align}
\sum_{\widetilde a \in A(\mathcal H^0)} \frac{G(0,0)}{\varphi(P)}
&= \frac{|A(\mathcal H^0)|}{P} \prod_{p \leq V} \left(1 -
\frac1{p}\right)^{-|\mathcal H^0|} \bar{\mathfrak S}_V(\mathcal
H^0)
\label{eq:12.17}\\
&= \prod_{p \leq V} \left(1 - \frac{\nu_p(\mathcal H^0)}{p}
\right) \left(1 - \frac1{p}\right)^{-|\mathcal H^0|} \cdot
\bar{\mathfrak S}_V (\mathcal H^0) = \mathfrak S(\mathcal H^0).
\nonumber
\end{align}
Inserting this into \eqref{eq:12.16} we obtain by \eqref{eq:12.3}
\begin{align}
&\widetilde{\widetilde S}_R (N; \mathcal H_1, \mathcal H_2,
\ell_1, \ell_2, P, h_0) :=
\sum_{\widetilde a \in A(\mathcal H)} \widetilde S_R(N;
\mathcal H_1, \mathcal H_2, \ell_1, \ell_2, P, \widetilde a, h_0)
\label{eq:12.18}\\
&= N \frac{{v + u \choose u} (\log R)^{d + v + u} \mathfrak
S(\mathcal H^0)}{(d + v + u)!} \left(1 + O \left(\frac{K \bar
d^* \log_2 N}{\log R} \right)\right)\nonumber \\
&\qquad + O \left(N \exp\left( - c\min \left(\sqrt{\log R},
\frac{\sqrt{\log N}}{\log_2 N}\right)\right)\right).
\nonumber
\end{align}

Now, from a brief examination of the values of the parameters $a,b,d$
in Cases 1, 2, 3 (after \eqref{eq:12.8}) and \eqref{eq:12.15},
we see that
\begin{equation}
\frac{{v + u \choose u} (\log R)^{d + v + u}}{(d + v + u)!} =
C_R(\ell_1, \ell_2, \mathcal H_1, \mathcal H_2, h_0) \cdot
\frac{{\ell_1 + \ell_2\choose \ell_1}
(\log R)^{r + \ell_1 + \ell_2}}{(r + \ell_1 + \ell_2)!}.
\label{eq:12.19}
\end{equation}
The relations \eqref{eq:12.18}--\eqref{eq:12.19} prove
Theorem~\ref{ujth:2}.

\section{The sum of the singular series $\mathfrak S(\mathcal H)$}
\label{sec:uj13}

Let
\begin{equation}
B_{\mathcal A}(k) = B(k) = \sum_{|\mathcal H| = k, \mathcal H
\subset \mathcal A} \mathfrak S(\mathcal H),
\label{8.1}
\end{equation}
where all sets
$\mathcal H = \{h_1, h_2, \dots, h_k\} \subseteq \mathcal A
\subseteq [1,N]$
are counted with $k!$ multiplicity according to
all possible permutations of $h_i$, and $|\mathcal A| = h$.

By Gallagher's theorem \cite{Ga} we have for fixed $k$
and $\mathcal A = [1,h]$
as $h \to \infty$
\begin{equation}
B_{\mathcal A}(k) = h^k \big(1 + O_{k,\varepsilon} (h^{-\frac12+\varepsilon})\big).
\label{8.2}
\end{equation}

This is not uniform in $k$ but up to some level $k \leq f(h)$
one could still show $B_{\mathcal A}(k) \sim h^k$.
However, we will use here a completely different approach.
We do not prove \eqref{8.2}, just (see Lemma~16) the weaker
relation that $B_{\mathcal A}(k)/h^k$ is, apart from a factor $1+o(1)$,
non-decreasing as a function of $k$, at least as long as $k =
o(h / \log h)$.
This result is fortunately completely sufficient for our purposes.

Further, our method is much more general and works for any set
$\mathcal A$ with $\mathcal A \subseteq [1,N]$, $|\mathcal A| = h$.

We remark that the asymptotic $B_{\mathcal A}(k) \sim h^k$ is
probably not true if $\mathcal A$ is arbitrary and even for
$\mathcal A = [1,h]$ it might fail if $k$ is as large as
$h/(\log h)^C$.

Let $c$ be an arbitrary small constant, $h, z, N$ and $Z$ sufficiently
large,
\begin{equation}
\gathered
k \leq \log N, \ \
h^2 \leq z = \log^5 N , \\
Z = P(z) = \prod_{p \leq z} p,\ \
Y = Y_z = \prod_{p \leq z} \left( 1 - \frac1 p \right)^{-1} \sim
e^\gamma \log z,
\endgathered
\label{8.3}
\end{equation}
\begin{equation}
Q := Q_z := \{n; (n, P(z)) = 1\}, \quad
M := \sum_{1 \leq n \leq Z, n \in Q} 1 = \frac Z Y.
\label{8.4}
\end{equation}
Then we have for a fixed set $\mathcal H$ consisting of $k$
distinct elements $h_i \in [1,N]$,
similarly to Section~6,  the density of $z$-quasi-prime
tuples of pattern~$\mathcal H$, using~(6.6):
\begin{equation}\begin{split}
R(\mathcal H) &:= \frac1 Z \sum^Z_{i = 1 \atop P_{\mathcal H}(i)
\in Q} 1\
= \ \prod_{p \leq z} \left(1 - \frac{\nu_p(\mathcal H)}{p}
\right)\
= \ Y^{-k} \prod_{p \leq z} \left(1 - \frac{\nu_p(\mathcal H)}{p}
\right) \left(1 - \frac1 p\right)^{-k} \\
& = Y^{-k} \mathfrak S(\mathcal H) \exp \left(O\left(k \sum_{p > z
\atop p \mid \Delta} \frac1p + k^2 \sum_{p > z \atop p \nmid
\Delta} \frac1{p^2}\right)\right) \\&
= Y^{-k} \mathfrak S(\mathcal H)\exp\left(O\left(k \sum_{p\mid
\Delta} \frac{\log p}{z \log z} + \frac{k^2}{z\log z}\right)\right)
\\
&= Y^{-k} \mathfrak S(\mathcal H) \exp \left(O\left(\frac{k^3
\log N}{z \log z}  + \frac{k^2}{z \log z}\right) \right) 
= Y^{-k} \mathfrak S(\mathcal H) \Big(1 + O\Big(\frac1{\log
N}\Big)\Big) , \label{8.5} \end{split}
\end{equation}
uniformly in $k, h, z, N $ satisfying \eqref{8.3}, if $c$ is fixed.
Let further
\begin{equation}
S^*_{\mathcal A}(k) := \frac1{h^k} \sum_{|\mathcal H| = k,
\mathcal H \subset \mathcal A} \mathfrak
S(\mathcal H) = \frac{B_{\mathcal A}(k)}{h^k}.
\label{8.6}
\end{equation}

\begin{lemma}
\label{lem:uj8}
If $k < \varepsilon(h) h/\log_2 N$, then
\begin{equation}
S^*_{\mathcal A}(k + 1) \geq S^*_{\mathcal A}(k)
\Big(1 + O \big(\varepsilon(h)\big) + O \Big(\frac1{\log N}\Big)\Big).
\label{8.7}
\end{equation}
\end{lemma}

\noindent
{\it Proof.}
For $i \in [1, Z]$ let
\begin{equation}
f_i = \sum_{j \atop i + a_j \in Q} 1, \
b_i = b_i(k) = f_i(f_i - 1) \dots (f_i - k+1).
\label{8.8}
\end{equation}
Then $b_i(k)$ is the number of all $k$-tuples of
$z$-quasiprimes of type $i + a_{j_\nu}$,
$a_{j_\nu} \in \mathcal A$
 ($\nu = 1,\dots, k$,
$1 \leq j_\nu \leq h$, $j_\nu$ distinct),
calculated with all $k!$
permutations, while $f_i$ is the number of
$z$-quasiprimes of the form $i + a_j$.
We have obviously for every pair $i,j \in [1, h]$
\begin{equation}
f_i \geq f_j \Leftrightarrow b_i \geq b_j,
\label{8.9}
\end{equation}
therefore
\begin{equation}
\frac1 Z \sum^Z_{i = 1} b_i f_i \geq \frac{\sum\limits^Z_{i=1}
f_i}{Z} \frac{\sum\limits^Z_{i = 1} b_i}{Z}.
\label{8.10}
\end{equation}
The above formula follows from
\begin{equation}
2 \left(Z\sum^Z_{i=1} b_i f_i - \sum^Z_{i=1} f_i \sum^Z_{i=1}
b_i \right) = \sum^Z_{i=1} \sum^Z_{j=1} (f_i - f_j) (b_i - b_j)
\geq 0.
\label{8.11}
\end{equation}
We have further $b_i(k + 1) = b_i(k)(f_i - k) = b_i f_i - k b_i$
and by calculating in two different ways how many times all pairs
$i,\mathcal H$ $(|\mathcal H| = k)$ satisfy the relation
$P_{\mathcal H}(i)  \in Q$ we
obtain
\begin{equation}
Z^{-1}\sum^Z_{i=1} b_{i}(k) =
Z^{-1}\sum^Z_{i=1} \sum_{\substack{|\mathcal H| = k
\\
P_{\mathcal H}(i) \in Q}} 1 = Z^{-1} \sum_{|\mathcal H| = k}
\sum^Z_{\substack{i=1\\
P_{\mathcal H}(i) \in Q}} 1 =
\sum_{|\mathcal H| = k} R(\mathcal H),
\label{8.12}
\end{equation}
while
\begin{equation}
\frac1 Z \sum^Z_{i=1} f_i = \frac{hM}{Z} = \frac h Y
.
\label{uj8.13}
\end{equation}

Thus \eqref{8.10} and \eqref{uj8.13} imply
by $b_i f_i = b_i(k+1) + kb_i$ that
\begin{equation}
\frac1 Z \sum^Z_{i=1} b_i(k + 1) + k \cdot \frac1 Z \sum^Z_{i=1}
b_i(k) \geq
\frac h Y  \cdot \frac1 Z
\sum^Z_{i=1} b_i(k).
\label{8.13}
\end{equation}
Hence, using \eqref{8.12}, we obtain
\begin{equation}
\sum_{|\mathcal H| = k+1} R(\mathcal H) \geq \left(\frac h Y -
k\right) \sum_{|\mathcal H| = k} R(\mathcal H).
\label{8.14}
\end{equation}
Multiplying by $Y^{k+1}$ on both sides, we obtain by \eqref{8.5}
\begin{equation}
\sum_{|\mathcal H| = k+1} \mathfrak S(\mathcal H) \geq h\left(1
+ O\left(\frac{kY}{h}\right)  + O \left(\frac1{\log N} \right)\right)
\sum_{|\mathcal H| = k} \mathfrak S(\mathcal H).
\label{8.15}
\end{equation}
Now dividing by $h^{k+1}$ on both sides we obtain \eqref{8.7} by
$Y \ll \log_2 N$.

\section{Proof of Theorem~\ref{th:A}}
\label{sec:13}

Theorems \ref{ujth:1} and \ref{ujth:2} allow us
 to express the quantity
$S'_R(N, K, \ell, P)$ in \eqref{eq:2.19}
in terms of
\begin{equation}
S^*_{\mathcal A}(k) = S^*(k) := \frac{B_{\mathcal A}(k)}{h^k}
:= \frac1{h^k} \sum_{|
\mathcal H| = k, \ \mathcal H \subset \mathcal A} \mathfrak S(\mathcal H),
\label{eq:13.1}
\end{equation}
where we consider two sets $\mathcal H$ and $\mathcal H'$
different if they contain the same elements in different
permutations.
The value of the parameter $k$ will be between $K$ and $2K + 1$
since in the application the sum \eqref{eq:13.1} will refer to
sums of type $\mathcal H = \mathcal H_1 \cup \mathcal H_2$,
$|\mathcal H_i| = K$, or to $\mathcal H^0$.

The derivation of the proof of Theorem~\ref{th:A} from our present
Theorems~\ref{ujth:1} and \ref{ujth:2} will be nearly
the same as that of the
main result (Theorem~3) of \cite{GPY} from
Propositions~1 and 2 in \cite{GPY},
which appears in Section~10 of \cite{GPY}, so we will be brief.
Although the restrictions for $K$ and $h$ will be quite
different here, nearly everything will be valid without any
change in the present case.
Our analysis refers now for the case $\nu = 1$ of Section~10 in
\cite{GPY}.

Let us choose, somewhat differently from \cite{GPY},
\begin{equation}
R = (3N)^{\Theta} = (3N)^{\frac14 - \xi}, \quad
\xi = c/\sqrt{\log N}, \quad V = \sqrt{\log N}
\label{eq:13.5}
\end{equation}
\begin{equation}
K = 16(\ell + 1)^2 = 16\varphi^{-2} \Longleftrightarrow \ell + 1
= \varphi^{-1} = \sqrt{K} / 4
\label{eq:13.6}
\end{equation}
\begin{equation}
x = \frac{K}{100} = \frac{\log R}{h} \Longleftrightarrow h =
\frac{100 \log R}{K} \left(\sim \frac{25\log N}{K}\right),
\label{eq:13.7}
\end{equation}
\begin{equation}
r_0 = (1 - 2\varphi) K , \quad
r_1 = (1 - \varphi)K
\label{eq:13.8}
\end{equation}
\begin{equation}
f(r) = {K\choose r}^2 \frac{x^r}{(r + 1) \dots(r + 2\ell)}, \quad
\bar r^* = \max(\sqrt K, K - r), \quad t(r) = \frac{\bar
r^*}{\varphi K}
\label{eq:13.9}
\end{equation}
and suppose that our crucial parameter $K$ satisfies
\begin{equation}
K \leq c_0 \frac{\sqrt{\log N}}{\log^2_2 N},
\label{eq:13.10}
\end{equation}
with a sufficiently small (explicitly calculable) absolute
constant $c_0$.

In the course of the proof of our present Theorem~\ref{th:A} (similar to
Section~10 of \cite{GPY}) a very important role is played by the fact
that although the sums evaluated in Theorems 4 and 5  depend on the actual choice of $\mathcal
H_1$ and $\mathcal H_2$, the asymptotic formulas for them
 depends just on the set $\mathcal H = \mathcal H_1 \cup
\mathcal H_2$ and on the size of $\mathcal H_1 \cap \mathcal H_2$.
On the other hand the size of the error terms may depend on the actual choice of
$\mathcal H_1, \mathcal H_2$ and $\mathcal H$.
This dependence is
made explicit in our present refined version, at
least in the sense that we show an asymptotic which is more
precise if $r = |\mathcal H_1
\cap \mathcal H_2|$ is near  $K = |\mathcal H_i|$.

We have seen in \cite{GPY} that taking any given set $\mathcal
H$
of given size $k = 2K - r \in [K, 2K]$,
we can write it in
\begin{equation}
(2K - r)!{K \choose r}^2 r!
\label{eq:13.11}
\end{equation}
ways as the union of two sets $\mathcal H_1$ and $\mathcal H_2$
of size $K$,
$|\mathcal H_1 \cap \mathcal H_2| = r$ if we consider sets
$\mathcal H_i$ and $\mathcal H'_i$ different when the permutation
of the same elements is different (cf.\ (10.4) of \cite{GPY}).
Now we can apply Theorems~\ref{ujth:1} and \ref{ujth:2}
in order to obtain similarly to
 Section~10 of \cite{GPY}
\begin{equation}
S'_R(N, K, \ell, P) = {2\ell\choose \ell} (\log R)^{2\ell}
P^*_{K, \ell}(x)
\label{eq:13.12}
\end{equation}
with
\begin{equation}
P^*_{K, \ell}(x) \geq \sum^K_{r = 0} f(r) S^*(2K - r)
\left(1 + O(\eta_2) + x \left( \frac{4K\left(1 -
\frac{\varphi}{2}\right)}{r + 2\ell + 1} - \frac1{\Theta} +
O(\eta_1)\right)\!\right)
\label{eq:13.13}
\end{equation}
where the error terms arising from Theorems~\ref{ujth:1},
\ref{ujth:2} and Lemma~\ref{lem:uj8}
are now
\begin{equation}
\eta_1 = \frac{K \bar r^*\log_2 N}{\log N} = \frac{4K^{3/2} t(r)
\log_2 N}{\log N}, \quad
\eta_2 = \frac1{\log^3_2 N},
\label{eq:13.14}
\end{equation}
and by our choice of $\Theta$ in \eqref{eq:13.5} we have
\begin{equation}
\frac1{\Theta} = 4 + O(\eta_3), \quad
\eta_3 = \frac1{\sqrt{\log N}}.
\label{eq:13.15}
\end{equation}

We will examine the quantity in parenthesis after $x$ in
\eqref{eq:13.13} which is clearly monotonic in $r$ (apart from
the error terms).
We have now
by \eqref{eq:13.6}--\eqref{eq:13.9}
for $r \leq r_1 = K - \varphi K = K - 4\sqrt K
\Leftrightarrow t(r) \geq 1$:
\begin{equation}
r + 2\ell + 1 < K - t(r) \varphi K + \frac{\sqrt K}{2} =
K\left(1 - t(r) \varphi + \frac{\varphi}{8}\right)
\label{eq:13.16}
\end{equation}
and therefore
\begin{align}
&\frac{4\left(1 - \frac{\varphi}{2}\right) K}{r + 2\ell + 1} -
\frac1{\Theta} + O (\eta_1) > 4\left(t(r) - \frac58\right)
\varphi + O (\eta_1 + \eta_3)
\label{eq:13.17}\\
&> \frac{16}{\sqrt K} \cdot \frac38 t(r) - \frac{C t(r)K^{3/2}
\log_2 N}{\log N} - \frac{C}{\sqrt{\log N}} > 0 \nonumber
\end{align}
by $t(r) \geq 1$ and \eqref{eq:13.10}.

On the other hand, as in (10.24)--(10.25)
of \cite{GPY}, the contribution of all terms
$r_2 > r_1$ to $P^*_{K, \ell}(x)$ is bounded by
\begin{equation}
e^{-\sqrt K} f(r_0) \max_{r > r_1} S^*(2K - r),
\label{eq:13.18}
\end{equation}
because $f(r)$ quickly decreases for $r > r_1$. (These
are the terms where the quantity in parenthesis after $x$ in \eqref{eq:13.13} may be negative.)
We have, for any $r > r_2$,
\begin{align}
\frac{f(r_2)}{f(r_0)}
&= \prod_{r_0 < r \leq r_2} \left( \frac{K - r + 1}{r} \cdot
\frac{\sqrt K}{10} \right)^2 \leq
\left(\frac{2\varphi K}{K - 2\varphi K} \cdot
\frac{\sqrt{K}}{10} \right)^{2\varphi K}
\label{mar15.16}\\
&\leq (0.81)^{8\sqrt K} = e^{-1.6 \sqrt{K}}. \nonumber
\end{align}
However, all terms $r \leq r_1$ have a positive contribution and
that of $r = r_0$ is at least
\begin{equation}
f(r_0) S^*(2K - r_0) \left(1 + O \left(\frac1{\log^3_2 N}\right)\right).
\label{eq:13.19}
\end{equation}
Now the quasi-monotonic property, Lemma~\ref{lem:uj8},
implies
\begin{equation}
\frac{S^*(2K - r_0)}{\max\limits_{r > r_1} S^*(2K - r)}
> e^{-(K - r_0) C/\log^3_2 N} = e^{-\frac{8C\sqrt K}{\log^3_2 N}}.
\label{eq:13.20}
\end{equation}
Consequently the positive term belonging to $r_0$ dominates all possibly
negative terms belonging to $r > r_1$ and therefore we have
\begin{equation}
P^*_{K, \ell} (x) > 0 \Longleftrightarrow S'_R(N, K, \ell, P) > 0.
\label{eq:13.21}
\end{equation}
This, by \eqref{eq:2.19}, proves the existence of some $n \in [N + 1, 2N]$ with
\begin{equation}
\sum_{p = n + a_\nu, a_\nu \in \mathcal A} \log p > \log (3N),
\label{eq:13.22}
\end{equation}
and thereby the existence of two primes $p', p'' \in [N + 1, 3N]$ with
\begin{equation}
0 \neq p'' - p' \in \mathcal A - \mathcal A.
\label{eq:13.23}
\end{equation}
This proves Theorem~\ref{th:A} if we choose $K$ maximal,
satisfying the restriction \eqref{eq:13.10}.

\renewcommand{\baselinestretch}{1}
\footnotesize
D. A. Goldston

Department of Mathematics

San Jose
State University

San Jose, CA 95192, USA

e-mail: goldston@math.sjsu.edu\\

J. Pintz

R\'enyi Mathematical Institute of the Hungarian Academy
of Sciences

Budapest, Re\'altanoda u. 13-15

H-1053

Hungary

e-mail: pintz@renyi.hu\\

C. Y. Y{\i}ld{\i}r{\i}m

Department of Mathematics

Bo\~{g}azi\c{c}i
University

Bebek

Istanbul 34342

Turkey

\&
Feza G\"ursey Enstit\"us\"u

\c Cengelk\"oy

Istanbul, P.K. 6,

81220 Turkey

e-mail: yalciny@boun.edu.tr

\end{document}